\documentclass[3p,times]{elsarticle}

\usepackage{ecrc}

\volume{00}

\firstpage{1}

\journalname{Computers $\&$ Mathematics with Applications}

\runauth{}

\jid{procs}

\jnltitlelogo{Computers $\&$ Mathematics with Applications}

\CopyrightLine{2011}{Published by Elsevier Ltd.}

\usepackage{amssymb}

\usepackage[figuresright]{rotating}

\usepackage{amsmath}
\usepackage{amsfonts}
\usepackage{cases}
\usepackage{float} 
\usepackage{subfig}
\usepackage[ruled]{algorithm2e} 
\newtheorem{remark}{Remark}

\usepackage{url}

\usepackage{changes}

\begin{document}

\begin{frontmatter}




\title{Neural Networks Based on Power Method and Inverse Power Method for Solving Linear Eigenvalue Problems}


\author[1]{Qihong Yang}
\ead{yangqh@stu.scu.edu.cn}

\author[1]{Yangtao Deng}
\ead{ytdeng1998@foxmail.com}

\author[1]{Yu Yang}
\ead{yuyang123@stu.scu.edu.cn}

\author[1]{Qiaolin He \corref{cor}}
\ead{qlhejenny@scu.edu.cn}
\author[1]{Shiquan Zhang \corref{cor}}
\ead{shiquanzhang@scu.edu.cn}

\cortext[cor]{Corresponding author}

\address[1]{School of Mathematics, Sichuan University, Chengdu, China}

\begin{abstract}
    In this article, we propose  two kinds of neural networks inspired by power method and inverse power method to solve linear eigenvalue problems. These neural networks share similar ideas with traditional methods, in which the differential operator is realized by automatic differentiation. The eigenfunction of the eigenvalue problem is learned by the neural network and the iterative algorithms are implemented by optimizing the specially defined loss function. The largest positive eigenvalue, smallest eigenvalue and interior eigenvalues with the given prior knowledge can be solved efficiently. We examine the applicability and accuracy of our methods in the numerical experiments in one dimension, two dimensions and higher dimensions. Numerical results show that accurate eigenvalue and eigenfunction approximations can be obtained  by our methods.
\end{abstract}


\end{frontmatter}



\section{Introduction}
\label{sec:introduction}
Recently, neural network has achieved remarkable success in solving partial differential equations (PDEs). Especially, deep Ritz method (DRM) proposed by E and Yu \citep{yu2018deep} and Physcis-informed neural network (PINN) proposed by Raissi et al.\citep{PINN} have attracted widespread attention. With the unprecedented availability of computational power, extensive research has been done on developing neural network for solving PDEs which can be traced back to 1990s \citep{lagaris1998artificial, lee1990neural, van1995neural}. These days, neural network has led to some remarkable results for solving a variety of problems, including heat transfer problems \citep{cai2021physics}, finance \citep{bai2022application}, uncertainty quantification \citep{gao2022wasserstein, oszkinat2022uncertainty, yang2019adversarial}, inverse problems \citep{chen2020physics, kadeethum2020physics} and so on.

How to effectively solve the eigenvalue problem is a very important issue, which has a \textcolor{blue}{wide} application in nuclear reactor physics field \cite{elhareef2022physics, buchan2013pod}, elastic acoustic problem \cite{diao2022spectral}, elastic viscoelastic composite structures \cite{chen2000integral}, etc. Traditionally, power method \citep{Golub1996} is used to find the dominant eigenvalue and eigenvector of a discretized eigenvalue equation, which depends on the mesh subdivision. When the degree of freedom is large, the method becomes less efficient. The development of neural network provides another way to solve the eigenvalue problem without mesh generation. 
To the best of the authors' knowledge, DRM \citep{yu2018deep} and PINN \citep{ben2020solving, ben2023deep} can solve both linear and nonlinear eigenvalue problems. However, DRM  is developed based on variational principle. In addition, DRM and PINN  do not work for all problems. In particular, they do not work for the eigenvalue problems with $\lambda \neq 0$ in higher dimension. As the first step, the main goal of the present article is to discuss the neural network method for solving linear eigenvalue problems. How to solve the nonlinear eigenvalue problems will be our future work. In nuclear engineering, the fundamental mode solution of K-eigenvalue problem \cite{bookalain2020} based on steady-state multi-group neutron diffusion theory is crucially required. Researchers have to numerically study the fundamental mode eigenvalue (called Keff) and the corresponding eigenvector for given geometry/material configurations. In recent work, PINN is also applied to solve the neutron diffusion equations \cite{WUzeyun2022, elhareef2022physics, JAGTAP2020113028, wang2022surrogate, YANG2023109656}. 
We see this work as a first step to develop the neural network to solve discontinuous interface K-eigenvalue problem in reactor physics.

As we known, there is a little research involved in solving the eigenvalue problems using neural network.
The Rayleigh quotient is used to solve the smallest eigenvalue problems in DRM \citep{yu2018deep}. In this article, the smallest eigenvalue is the minimum modulus eigenvalue, which is defined in Section 2.1. The authors transformed the original eigenvalue problem to the function which is known as the Rayleigh quotient using the variational principle for the smallest eigenvalue. The main idea of the Rayleigh quotient is the fact that it gives the range of eigenvalues of the operator. Then, the authors minimized the function and got the smallest eigenvalue and the associated eigenfunction which is expressed by the neural network. Following the similar idea, there are some works \citep{ben2020solving, ben2023deep} using the Rayleigh quotient to construct the function without variation, which is obtained by PINN. 

In \citep{jin2022physics}, several regularization terms are added into the loss function, the neural network will try to learn the smallest eigenvalue by minimizing the loss function. In \citep{han2020solving}, the eigenvalue problem is reformulated as a fixed point problem of the semigroup flow induced by the operator, whose solution can be represented by Feynman-Kac formula in terms of forward-backward stochastic differential equations (FBSDEs), where diffusion Monte Carlo method is used and the eigenfunction is approximated through neural network ansatz. 
Unfortunately, 
the methods of adding regularization terms into the loss function are difficult to be used in practical applications, since it is difficult to find the convergence point in the learning curve of a neural network. To some extent, this approach relies on the choice of hyper-parameters.
 For instance, the regularization term $e^{-\lambda+c}$ relies heavily on the constant $c$ in \citep{jin2022physics}, where $\lambda$ is the eigenvalue.

In addition, all methods discussed above are 
 trying to  minimize a loss function which may represent the eigenvalue to some extent to solve the  eigenvalue problems, where the neural network may fail to learn the smallest eigenvalue in some cases.
 The formulation of the loss function includes many terms so that it is  difficult to balance the weights of them for the neural network. Also, little research related to neural network studies on solving for the largest positive eigenvalue and the interior eigenvalues. 

Therefore, we are going to develop new methods to overcome the shortcomings of directly minimizing eigenvalue. It is well known that the power method is widely used to find the maximum eigenvalue of $n$ by $n$ matrix and the corresponding eigenvector. 
Our motivation is to take advantage of the power method and the neural network to explore efficient algorithms for eigenvalue problems.
In this work, we propose the power method neural network (PMNN) and inverse power method neural network (IPMNN) to solve eigenvalue problems with the largest positive eigenvalue, the smallest eigenvalue and the interior eigenvalues when an approximation is given.
In similar spirits of traditional methods, our methods focus on linear differential operators and iteratively approximate the exact eigenvalue and eigenfunction. 
Instead of discretizing the differential operator to a matrix system, we use automatic differentiation (AD) \citep{baydin2018automatic} to represent the operator. We use the neural network to learn the eigenfunction through optimizing the specially defined loss function,  which can be evaluated anywhere in the computational domain. This is a big difference from the traditional method. 

This article is organized as follows. In Section \ref{sec:preliminaries}, the power method and the inverse power method are reviewed. In Section \ref{sec:methods}, we propose two new methods PMNN and IPMNN . In Section \ref{sec:experiments}, numerical experiments are presented to verify our methods. 
Finally, the conclusions and future work are given in Section \ref{sec:conclusions}.

\section{Preliminaries}
\label{sec:preliminaries}

\subsection{Eigenvalue Problems}
\label{sec:eigenvalue_problems}
In this work, we focus on self-adjoint operators $\mathcal{L}$. Suppose that we have the following PDE,

       \begin{numcases}{}
          \mathcal{L}u= \lambda u, \  \mbox{in} \enspace \Omega,      \label{eq:eigenvalue_equation1}    \\
               \mathcal{B}u = g,       \  \mbox{on} \enspace \partial \Omega,  \label{eq:eigenvalue_equation2}
       \end{numcases}
where the domain $\Omega \subset \mathbb{R}^d$, $\mathcal{L}$ and $\mathcal{B}$ are two differential operators acting on the functions defined in the interior of $\Omega$ and on the boundary $\partial \Omega$, respectively. The $(u, \lambda)$ is an eigenpair of the operator $\mathcal{L}$,  
where $u$ is the eigenfunction of $\mathcal{L}$ and $\lambda$ is the corresponding eigenvalue. 
Equation \eqref{eq:eigenvalue_equation1}-\eqref{eq:eigenvalue_equation2} presents a generic form of the linear eigenvalue problem. Equation \eqref{eq:eigenvalue_equation1} is utilized to construct the Rayleigh quotient as shown below. The boundary condition \eqref{eq:eigenvalue_equation2} is the constraint necessary for the eigenvalue problem. For simplicity, we only consider the Dirichlet boundary condition and the periodic boundary condition in this article.

Traditionally, finite difference method (FDM) \citep{truhlar1972finite, simos1997finite}, finite element method \cite{ishihara1977convergence, ishihara1978mixed}, finite volume method \cite{liang2001finite, dai2011finite} and spectral method \citep{talbot2005application, atkinson2009spectral} have been fully developed to solve equation \eqref{eq:eigenvalue_equation1}-\eqref{eq:eigenvalue_equation2}. The matrix equation is formulated by the above methods, which can be solved by the power method or 
the inverse power method to obtain the dominant eigenvalue or smallest eigenvalue. In the following section, we will briefly introduce the power method and the inverse power method \citep{Golub1996,DW2002}.

\begin{remark}
    If the differential operator $\mathcal{L}$ is singular which means that the smallest eigenvalue is $\lambda = 0$, it is better to shift the operator $\mathcal{L}$ as shown in \eqref{eq:Lu_equation1} and solve the new problem using  Algorithm 4.
\end{remark}

\subsection{Power Method}
\label{sec:power_method}
The power method is widely used to find the largest eigenvalue (in absolute value) and the corresponding eigenvector of the discretized matrix. Such an eigenvalue is also called a dominant eigenvalue. Suppose $A \in \mathbb{R}^{n \times n}$ is a diagonalizable $n \times n$ matrix, which has $n$ eigenvalues $\lambda_1, \lambda_2, \cdots, \lambda_n$, with $\lambda_1$ being the dominant eigenvalue. The eigenvalues can be listed as $\lvert \lambda_1 \rvert > \lvert \lambda_2 \rvert \geq \lvert \lambda_3 \rvert \geq \cdots \geq \lvert \lambda_n \rvert \geq 0$. 
Let $V \in \mathbb{R}^{n \times n}$ be the matrix composed of its linearly independent unit eigenvector $\boldsymbol{v_i}$ and the corresponding  eigenvalue is $\lambda_i$, for $i=1, 2, \cdots, n$. There 
exists some constants $c_1, c_2, \cdots, c_n$ such that an arbitrary vector $\boldsymbol{v} \in \mathbb{R}^{n}$ can be uniquely expressed as
\begin{equation}
    \label{eq:x_vector}
    \begin{array}{r@{}l}
        \begin{aligned}
             & \boldsymbol{v}=c_1 \boldsymbol{v_1} + c_2 \boldsymbol{v_2} + \cdots + c_n \boldsymbol{v_n}. \\
        \end{aligned}
    \end{array}
\end{equation}
 We assume that $c_1 \neq 0$ and multiply \eqref{eq:x_vector} by $A$. Since $A \boldsymbol{v_i} = \lambda_i \boldsymbol{v_i}$,  we have 
\begin{equation}
    \label{eq:x_mul_A1}
    \begin{array}{r@{}l}
        \begin{aligned}
             & A \boldsymbol{v}=c_1 \lambda_1 \boldsymbol{v_1} + c_2 \lambda_2 \boldsymbol{v_2} + \cdots + c_n \lambda_n \boldsymbol{v_n}. 
        \end{aligned}
    \end{array}
\end{equation}
It is easy to obtain that
\begin{equation}
    \label{eq:x_mul_A2}
    \begin{array}{r@{}l}
        \begin{aligned}
            A^k \boldsymbol{v} & = c_1 \lambda_1^k \boldsymbol{v_1} + c_2 \lambda_2^k \boldsymbol{v_2} + \cdots + c_n^k \lambda_n \boldsymbol{v_n}                                                                   \\
                               & = c_1 \lambda_1 ^k \left(\boldsymbol{v_1} + \frac{c_2}{c_1} \frac{\lambda_2^k}{\lambda_1^k} \boldsymbol{v_2} + \cdots  + \frac{c_n}{c_1} \frac{\lambda_n^k}{\lambda_1^k} \boldsymbol{v_n}\right).
        \end{aligned}
    \end{array}
\end{equation}

Since $\lambda_1$ is the dominant eigenvalue, 
the ratios $\lvert \lambda_2\lvert / \lvert \lambda_1 \lvert, \lvert\lambda_3 / \lvert  \lambda_1 \lvert, \cdots, \lvert \lambda_n \lvert / \lvert \lambda_1 \lvert$ are all strictly less than 1. Hence for a sufficiently large $k$, $\lvert \lambda_1 \lvert^k$ is significantly larger than $\lvert \lambda_i \lvert^k (2 \leq i \leq n)$, then the ratio $\lvert \lambda_i \lvert^k/ \lvert \lambda_1 \lvert^k$ approaches 0 as $k \to \infty$. Thus, the term $c_1 \lambda_1^k \boldsymbol{v_1}$ dominates the expression for $A^k \boldsymbol{v}$ for large value of $k$.  We normalize $A^k \boldsymbol{v}$ and obtain $\boldsymbol{u}=(A^k \boldsymbol{v})/\lVert A^k \boldsymbol{v} \rVert \approx (c_1 \lambda_1^k \boldsymbol{v_1})/ \lVert c_1 \lambda_1^k \boldsymbol{v_1} \rVert$, which is a scalar mutilple of $\boldsymbol{v_1}$. Thus, $\boldsymbol{u}$ is a unit eigenvector corresponding to the dominant eigenvalue $\lambda_1$.
Therefore, we have $A\boldsymbol{u} \approx \lambda_1 \boldsymbol{u}$. 
The sign of $\lambda_1$ is determined by checking whether $A\boldsymbol{u}$ is in the same direction as $\boldsymbol{u}$ or not. Actually, the approximated eigenvalue in practice  is usually obtained via Rayleigh Quotient $\boldsymbol{u}^T A \boldsymbol{u}$ instead of $\lVert A\boldsymbol{u} \rVert$. The implementation details of power method is shown in Algorithm \ref{algo:power_method}.

\begin{algorithm}[htp]
    \caption{Power method for finding the dominant eigenvalue of  the matrix $A$}
    \label{algo:power_method}
    Let $A \in \mathbb{R}^{n \times n}$  be an $n \times n$ matrix. \\
    \textbf{Step 1:} Choose an arbitrary unit vector 
     $\boldsymbol{u_0} \in \mathbb{R}^{n}$, the maximum number of iterations $k_{max}$ and the stopping criterion $\epsilon$. \\
    \textbf{Step 2:} \\
    \For{$k=1, 2, \cdots, k_{max}$}{
        $\boldsymbol{p_k} = A \boldsymbol{u_{k-1}}$. \\
        $\boldsymbol{u_k} = \boldsymbol{p_k}/ \lVert \boldsymbol{p_k} \rVert$. \\
       \uIf{$\lVert \boldsymbol{u_k} - \boldsymbol{u_{k-1}} \rVert < \epsilon \space$}{
           The stopping criterion is met, the iteration is stopped. \\
        }\Else{
            Let $\boldsymbol{u_{k-1}} = \boldsymbol{u_k}$. \\
        }
    }
     The dominant eigenvalue $\lambda$ is obtained by Rayleigh quotient, \\
     $\lambda = \boldsymbol{u_k}^T A \boldsymbol{u_k}$.
\end{algorithm}

\subsection{Inverse Power Method}
\label{sec:inverse_power_method}
The inverse power method is widely used to find the smallest eigenvalue (in absolute value) and the corresponding eigenvector. 
Suppose $\lambda_1$ being the smallest eigenvalue and the eigenvalues can be listed as $0 < \lvert \lambda_1 \rvert < \lvert \lambda_2 \rvert \leq \lvert \lambda_3 \rvert \leq \cdots \leq \lvert \lambda_n \rvert$. Therefore, $A$ is a nonsingular matrix and  the inverse of it $A^{-1}$ exists.
Let $V \in \mathbb{R}^{n \times n}$ be the matrix composed of its linearly independent eigenvectors $\boldsymbol{v_i}$ and the corresponding  eigenvalue is $\lambda_i$, for $i=1, 2, \cdots, n$.

The same as power method, an arbitrary vector $\boldsymbol{v} \in \mathbb{R}^{n}$ can be uniquely expressed as \eqref{eq:x_vector}. 
It is easy to see that $\lvert 1/\lambda_1 \rvert > \lvert 1/\lambda_2 \rvert \geq \lvert 1/\lambda_3 \rvert \geq \cdots \geq \lvert 1/\lambda_n \rvert > 0$ and $1/\lambda_1, 1/\lambda_2, \cdots, 1/\lambda_n$ are the eigenvalues of the matrix $A^{-1}$. Equivalently, $A^{-1}\boldsymbol{u}=\lambda \boldsymbol{u}$ will be solved. 
Hence, the problem is to solve for the dominant eigenvalue of the matrix $A^{-1}$. The implementation details of inverse power method is shown in Algorithm \ref{algo:inverse_power_method}.

\begin{algorithm}[htp]
    \caption{Inverse power method for finding the smallest eigenvalue of the matrix $A$}
    \label{algo:inverse_power_method}
    Let $A \in \mathbb{R}^{n \times n}$ be an $n \times n$ matrix. \\
    \textbf{Step 1:} Choose an arbitrary unit vector $\boldsymbol{u_0} \in \mathbb{R}^{n}$, the maximum number of iterations $k_{max}$ and the stopping criterion $\epsilon$. \\
    \textbf{Step 2:} \\
    \For{$k=1, 2, \cdots, k_{max}$}{
       $\boldsymbol{p_k} = A^{-1} \boldsymbol{u_{k-1}}$ (LU decomposition is used). \\
       $\boldsymbol{u_k} = \boldsymbol{p_k}/ \lVert \boldsymbol{p_k} \rVert$.  \\
        \uIf{$\lVert \boldsymbol{u_k} - \boldsymbol{u_{k-1}} \rVert < \epsilon \space$}{
            The stopping criterion is met, the iteration can be stopped. \\
        }\Else{
            Let $\boldsymbol{u_{k-1}} = \boldsymbol{u_k}$. \\
        }
    }
    The smallest eigenvalue $\lambda$ is obtained by Rayleigh quotient, \\
    $\lambda = \boldsymbol{u_k}^T A \boldsymbol{u_k}$. 
\end{algorithm}

\section{Methodologies}
\label{sec:methods}
In this work, we aim to solve the eigenvalue problems using neural network.
We propose two architectures to solve the eigenvalue problems, called PMNN and IPMNN.
Especially, considering the main idea of the shifted inverse
power method, IPMNN can be used to obtain the interior eigenvalues when some prior knowledge $\alpha$ is given, which is attained by shifting the differential operator.

\subsection{Neural Networks for Solving Eigenvalue Problems}
We use $\mathcal{N}^{\theta}(\boldsymbol{x})$ to denote a neural network and the eigenfunction $u(\boldsymbol{x})$ can be represented by $\mathcal{U}(\boldsymbol{x})=\mathcal{N}^{\theta}(\boldsymbol{x})$, where $\theta$ denotes the parameters of the neural network.
For equations \eqref{eq:eigenvalue_equation1}--\eqref{eq:eigenvalue_equation2}, we propose two methods to solve for the largest positive eigenvalue and the smallest eigenvalue. In our methods, the function $u$ is expressed by neural network,
where boundary conditions will be discussed in Section \ref{sec:enforcement_boundary_conditions}.
The eigenfunction  expressed by neural network can predict value for any points in $\Omega$.

\subsubsection{Power Method Neural Network}
\label{sec:pmnn}
Inspired by the idea of power method, we propose PMNN to solve for the largest positive eigenvalue and the associated eigenfunction.
In PMNN,  we use the neural network $\mathcal{N}^{\theta}$ to represent the approximated eigenfunction $\mathcal{U}_{k-1}$ in the $k$-th iteration and calculate $\mathcal{P}_{k}$ by equation \eqref{eq:key_equation_pmnn1}, which is an analogue to $\boldsymbol{p_k} = A\boldsymbol{u_{k-1}}$ in Algorithm \ref{algo:power_method}. Similar as normalization $\boldsymbol{u_k} = \boldsymbol{p_k}/ \lVert \boldsymbol{p_k} \rVert$ in Algorithm \ref{algo:power_method}, the $\mathcal{P}_{k}$ should also be normalized in \eqref{eq:key_equation_pmnn2}. 
Here, the $\mathcal{U}_{k}$ is the approximated eigenfunction in the $k$-th iteration which can be obtained by equation \eqref{eq:key_equation_pmnn2}.
Different from the original power method, $\mathcal{L}$ is an operator, which is realized by AD. Therefore, the term $\mathcal{L}\mathcal{U}_{k-1}$ can be computed by AD.

\begin{eqnarray}
            & & \mathcal{P}_{k} = \mathcal{L}\mathcal{U}_{k-1},\label{eq:key_equation_pmnn1}\\ 
            & & \mathcal{U}_{k} = \frac{\mathcal{P}_{k}}{\lVert \mathcal{P}_k \rVert}\label{eq:key_equation_pmnn2}.
\end{eqnarray}


Although $\mathcal{U}_{k}$ is obtained by \eqref{eq:key_equation_pmnn1}--\eqref{eq:key_equation_pmnn2}, it is much different from traditional power method and we can not  assign the values of $\mathcal{U}_{k}$ to the neural network. 
Therefore,  we define a loss function in equation \eqref{eq:loss_pmnn} to represent the convergence condition in the power method, 
where $\boldsymbol{x}_i \in S$, which is the data set attained by using random sampling algorithm, and $N$ denotes the number of points in the data set $S$.
Through the defined loss function, the eigenfunction will be approximated iteratively.
The main idea behind this equation is that we do not need the neural network to precisely compute $\mathcal{U}_k$ in equation \eqref{eq:key_equation_pmnn2} in the next step. The approximation is done step by step.
\begin{equation}
    \label{eq:loss_pmnn}
    \begin{array}{r@{}l}
        \begin{aligned}
            loss_{pmnn}(\theta) = \frac{1}{N}\sum_{i=1}^{N}[\mathcal{U}_{k-1}(\boldsymbol{x}_i)-\mathcal{U}_{k}(\boldsymbol{x}_i)]^2.
        \end{aligned}
    \end{array}
\end{equation}


When the neural network gets convergence, we obtain the largest positive eigenvalue and associated eigenfunction expressed by the neural network. For the process of this method, one can refer to Algorithm \ref{algo:pmnn}.

\begin{algorithm}[htp]
    \caption{PMNN for finding the largest positive eigenvalue}
    \label{algo:pmnn}
    Give $N$ the number of points for training neural network, $N_{epoch}$ the maximum number of epochs and the stopping criterion $\epsilon$. \\
    \textbf{Step 1:} Build data set $S$ for training using random sampling algorithm. \\
    \textbf{Step 2:} Initialize a neural network with random initialization of parameters.  \\
    \For{$k=1, 2, \cdots, N_{epoch}$}{
        Input all points in $S$ into neural network $\mathcal{N}^{\theta}$.\\
        Let $\mathcal{U}_{k-1}(\boldsymbol{x}_i) = \mathcal{N}^{\theta}(\boldsymbol{x}_i)$, where $\boldsymbol{x}_i \in S$. \\
        $\mathcal{P}_{k} = \mathcal{L}\mathcal{U}_{k-1}$ (using AD).\\
        $\mathcal{U}_{k} = \frac{\mathcal{P}_k}{\lVert \mathcal{P}_{k} \rVert}$. \\
        Calculate $loss_{pmnn}(\theta) = \frac{1}{N}\sum_{i=1}^{N}[\mathcal{U}_{k-1}(\boldsymbol{x}_i)-\mathcal{U}_{k}(\boldsymbol{x}_i)]^2$, where       $\boldsymbol{x}_i \in S$. \\
        Update parameters of neural network using optimizer. \\
        \If{$loss_{pmnn} < \epsilon$}{
            Record the  eigenvalue and eigenfunction, \\
            $\lambda=\frac{<\mathcal{L}\mathcal{U}_{k-1},\mathcal{U}_{k-1}>}{<\mathcal{U}_{k-1}, \mathcal{U}_{k-1}>}$. \\
            The stopping criterion is met, the iteration can be stopped. \\
        }
    }
\end{algorithm}

\begin{remark}
    Compared with the traditional power method, 
    we use the assignment expression $\mathcal{U}_{k-1}(\boldsymbol{x}_i) = \mathcal{N}^{\theta}(\boldsymbol{x}_i)$ without normalization, where $\boldsymbol{x}_i \in S$. The reason is that 
    $\mathcal{U}_{k-1}$ will be close to $\mathcal{U}_{k}$ after training the neural network, where $\mathcal{U}_{k}$ is a normalized expression in the this Algorithm. Obviously, $\mathcal{U}_{k-1}$ is also normalized during the training process.
\end{remark}

\begin{remark}
    In both the power method and PMNN, we can obtain $\mathbf{u}_k$ and $\mathcal{U}_k$ from $\mathbf{u}_{k-1}$ and $\mathcal{U}_{k-1}$. But the vector $\mathbf{u}_k$ is assigned to $\mathbf{u}_{k-1}$ in the next iteration in the power method and it is impossible to do so in the neural network. Therefore, we define a loss function to propel the neural network to learn $\mathcal{U}_k$, which provide a correct direction for optimization.
\end{remark}

\subsubsection{Inverse Power Method Neural Network}
\label{sec:ipmnn}
Inspired by the idea of inverse power method, we propose IPMNN to solve for the smallest eigenvalue and the associated eigenfunction.
In IPMNN, it is different from PMNN that we use the neural network $\mathcal{N}^{\theta}$ to represent the approximated eigenfunction $\mathcal{U}_{k-1}$ in the $k$-th iteration. Here, we use the neural network $\mathcal{N}^{\theta}$ to represent the approximated  eigenfunction $\mathcal{U}_{k}$ in the $k$-th iteration in equation \eqref{eq:key_equation_ipmnn}
 which is an analogue to equation \eqref{eq:key_equation_inverse_power_method} in Algorithm \ref{algo:inverse_power_method}. 
Given $\mathcal{U}_{k-1}$ which is from the last iterative step and following the main idea of inverse power method, we need to solve $\mathcal{U}_{k}$ by $\mathcal{P}_k=\mathcal{L}^{-1}\mathcal{U}_{k-1}$ and $\mathcal{U}_{k} = \mathcal{P}_{k} / \lVert \mathcal{P}_{k} \rVert$. However, it is difficult to get the inverse operator $\mathcal{L}^{-1}$ of the differential operator $\mathcal{L}$. Therefore, we are going to get $\mathcal{U}_{k}$ without knowing the inverse operator. 
 Since $\mathcal{L}$ is realized by AD in PMNN, the term $\mathcal{L}\mathcal{U}_k$ is computed by AD in equation \eqref{eq:key_equation_ipmnn} similarly.
\begin{equation}
    \label{eq:key_equation_inverse_power_method}
    \begin{array}{r@{}l}
        \left\{
        \begin{aligned}
            \boldsymbol{p_k} & = A^{-1}\boldsymbol{w_{k-1}},                              \\
            \boldsymbol{w_k} & = \frac{\boldsymbol{p_k}}{\lVert \boldsymbol{p_k} \rVert}. \\
        \end{aligned}
        \right.  \quad  (\mbox{in} \ \mbox{Algorithm} \ \ref{algo:inverse_power_method})
    \end{array}
\end{equation}
\begin{equation}
  \label{eq:key_equation_ipmnn}
   \begin{array}{r@{}l}
        \begin{aligned}
            \frac{\mathcal{L}\mathcal{U}_{k}}{\lVert \mathcal{L}\mathcal{U}_{k} \rVert} = \mathcal{U}_{k-1}. \\
        \end{aligned}
    \end{array}
\end{equation}

Since it is impossible to directly compute $\mathcal{P}_k=\mathcal{L}^{-1}\mathcal{U}_{k-1}$
and the neural network is used to express $\mathcal{U}_k$, the equality of Equation \eqref{eq:key_equation_ipmnn} will not be satisfied exactly. Therefore, 
we define a loss function in equation \eqref{eq:loss_ipmnn} to represent the convergence condition in the inverse power method and relieve this problem, 
where $\boldsymbol{x}_i \in S$, which is the data set attained by using random sampling algorithm, and $N$ denotes the number of points in the data set $S$. 
The main idea is that we do not need to calculate $\mathcal{L}^{-1}$ and the eigenfunction will be approximated iteratively through minimizing the defined loss \eqref{eq:loss_ipmnn} to approach the equation \eqref{eq:key_equation_ipmnn}.
\begin{equation}
    \label{eq:loss_ipmnn}
    \begin{array}{r@{}l}
        \begin{aligned}
            loss_{ipmnn}(\theta) = \frac{1}{N}\sum_{i=1}^{N}\left(\frac{\mathcal{L}\mathcal{U}_{k}(\boldsymbol{x}_i)}{\lVert \mathcal{L}\mathcal{U}_{k} \rVert}-\mathcal{U}_{k-1}(\boldsymbol{x}_i)\right)^2.
        \end{aligned}
    \end{array}
\end{equation}

When the neural network gets convergence, we  obtain the smallest eigenvalue and the associated eigenfunction expressed by the neural network. For the process of this method, one can refer to Algorithm \ref{algo:ipmnn}.

\begin{algorithm}[htp]
    \caption{IPMNN for finding the smallest eigenvalue}
    \label{algo:ipmnn}
    Give $N$ the number of points for training neural network, $N_{epoch}$ the maximum number of epochs and the stopping criterion $\epsilon$. \\
    \textbf{Step 1:} Build data set $S$ for training using random sampling algorithm. \\
    \textbf{Step 2:} Choose an arbitrary normalized function $u_0$, and let $\mathcal{U}_{0} = u_0$. \\
    \textbf{Step 3:}  Initialize a neural network with random initialization of parameters. \\
    \For{$k=1, 2, \cdots, N_{epoch}$}{
        Input all points in $S$ into neural network $\mathcal{N}^{\theta}$. \\
        Let $\mathcal{U}_k(\boldsymbol{x}_i)=\mathcal{N}^{\theta}(\boldsymbol{x}_i)$, where $\boldsymbol{x}_i \in S$.  \\
        Compute $\mathcal{L}\mathcal{U}_{k}$ using AD. \\
        Calculate $loss_{ipmnn}(\theta) = \frac{1}{N}\sum_{i=1}^{N}\left(\frac{\mathcal{L}\mathcal{U}_{k}(\boldsymbol{x}_i)}{\lVert \mathcal{L}\mathcal{U}_{k} \rVert}-\mathcal{U}_{k-1}(\boldsymbol{x}_i)\right)^2$,  where $\boldsymbol{x}_i \in S$. \\
        Update parameters of neural network using optimizer. \\
        $\mathcal{U}_{k-1} = \frac{\mathcal{U}_{k}}{\lVert \mathcal{U}_{k} \rVert}$. \\
        \If{$loss_{ipmnn} < \epsilon$}{
            Record the eigenvalue and eigenfunction, \\
            $\lambda=\frac{<\mathcal{L}\mathcal{U}_{k}, \mathcal{U}_{k}>}{<\mathcal{U}_{k}, \mathcal{U}_{k}>}$. \\
            The stopping criterion is met, the iteration can be stopped. \\
        }
    }
\end{algorithm}

If the operator $\mathcal{L}$ is shifted,  IPMNN can be used to obtain the interior eigenvalues when some prior knowledge $\alpha$ is given.
The main idea comes from the shifted inverse power method.  
We subtract $\alpha u$ from both sides of equation \eqref{eq:eigenvalue_equation1}. Then we can get equation \eqref{eq:Lu_equation1}, where $\widetilde{\mathcal{L}}=\mathcal{L}-\alpha \mathcal{I}$, $\widetilde{\lambda} =\lambda-\alpha$, and $\mathcal{I}$ is the identity operator. Therefore, the original problem can be rewritten as equation \eqref{eq:Lu_equation1}, which can be solved by IPMNN.
\begin{eqnarray}
               & &  \widetilde{\mathcal{L}} u  =   \widetilde{\lambda} u, \  \mbox{in} \enspace \Omega.  \label{eq:Lu_equation1}
\end{eqnarray}
So, we can obtain the eigenvalue which is close to $\alpha$ and the corresponding eigenfunction of problem \eqref{eq:eigenvalue_equation1}. Our numerical results of the harmonic eigenvalue problem verify the efficiency of IPMNN for different $\alpha$.

\begin{remark}
Although the eigenfunction is represented by the neural network in both PMNN and IPMNN, we use the neural network $\mathcal{N}^{\theta}$ to represent the approximated eigenfunction $\mathcal{U}_{k-1}$ in the $k$-th iteration in PMNN and to represent the approximated  eigenfunction $\mathcal{U}_{k}$ in the $k$-th iteration in IPMNN. Therefore, we do not need to give the initial eigenfunction $\mathcal{U}_0$ in PMNN and it is necessary to give the initial eigenfunction $\mathcal{U}_0$ in IPMNN.
\end{remark}



\subsection{Enforcement of Boundary Conditions}
\label{sec:enforcement_boundary_conditions}
As we know, the implementation of boundary conditions is very important for PDEs system  \citep{evans2010partial}. 
For simplicity, we only consider the Dirichlet boundary condition and the periodic boundary condition in this work.
It is easy to use the strategies in \citep{berg2018} and \citep{lyu2020enforcing} to enforce the exact boundary conditions. In our future work, more complicated domains or more complicated boundary conditions will be considered. We may enforce the boundary conditions by adding penalty terms and using data points on the boundary. 
For example, we may define a loss function $loss_b = \sum_{i=1}^{N_b} \lvert \mathcal{B} u(\boldsymbol{x}_i) - g (\boldsymbol{x}_i)\rvert ^2$, where $N_b$ is the number of sampling points on $\partial \Omega$ and $\boldsymbol{x}_i$ is the point in the sampling set $\{\boldsymbol{x}_i\}_{i=1}^{N_b}$.

\subsubsection{Enforcement of Dirichlet Boundary Condition}
\label{sec:enforcement_dirichlet_boundary_conditions}
For the Dirichlet boundary condition, we can directly access the values of solution on the boundary. The boundary condition \eqref{eq:eigenvalue_equation2} can be written as following,
\begin{equation}
    \label{eq:dirichelte_boundary_condition}
    \begin{array}{r@{}l}
        \begin{aligned}
             & u(\boldsymbol{x}) = g(\boldsymbol{x}), &  & \forall \boldsymbol{x} \in \partial \Omega.
        \end{aligned}
    \end{array}
\end{equation}
Therefore, it is straightforward to design a distance function $\phi(\boldsymbol{x})$, 
\begin{equation}
    \label{eq:phi_x}
    \begin{array}{r@{}l}
         & \phi(\boldsymbol{x}) \left\{
        \begin{aligned}
             & =0,     &  & \forall \boldsymbol{x} \in \partial \Omega, \\
             & \neq 0, &  & \forall \boldsymbol{x} \in \Omega.
        \end{aligned}
        \right.
    \end{array}
\end{equation}
Suppose $\mathcal{N}^\theta (\boldsymbol{x})$ be the output of neural networks, the solution $\mathcal{U} (\boldsymbol{x})$ is constructed as
\begin{equation}
    \label{eq:u_enfocing_dirichlet}
    \begin{array}{r@{}l}
        \begin{aligned}
             & \mathcal{U} (\boldsymbol{x}) = \phi(\boldsymbol{x}) \mathcal{N}^\theta(\boldsymbol{x}) + G(\boldsymbol{x}), &  & \forall \boldsymbol{x} \in \overline{\Omega},
        \end{aligned}
    \end{array}
\end{equation}
where $G(\boldsymbol{x})$ is a (smooth) extension of $g(\boldsymbol{x})$ in $\Omega$ and $G(\boldsymbol{x}) = g(\boldsymbol{x}), \  \forall \mathbf{x} \in \partial \Omega$. For instance, suppose $\Omega=(0, 1)$ and $g(x)=x^2$ when $x=0$ or $x=1$. Then, we can simply define $G(x)=x^2, \  \forall \mathbf{x} \in \overline{\Omega}$,  using the approach in \citep{lyu2020enforcing}.

\subsubsection{Enforcement of Periodic Boundary Condition}
\label{sec:enforcement_periodic_boundary_conditions}
\begin{figure}[htp]
    \centering
    \includegraphics[width=1\textwidth]{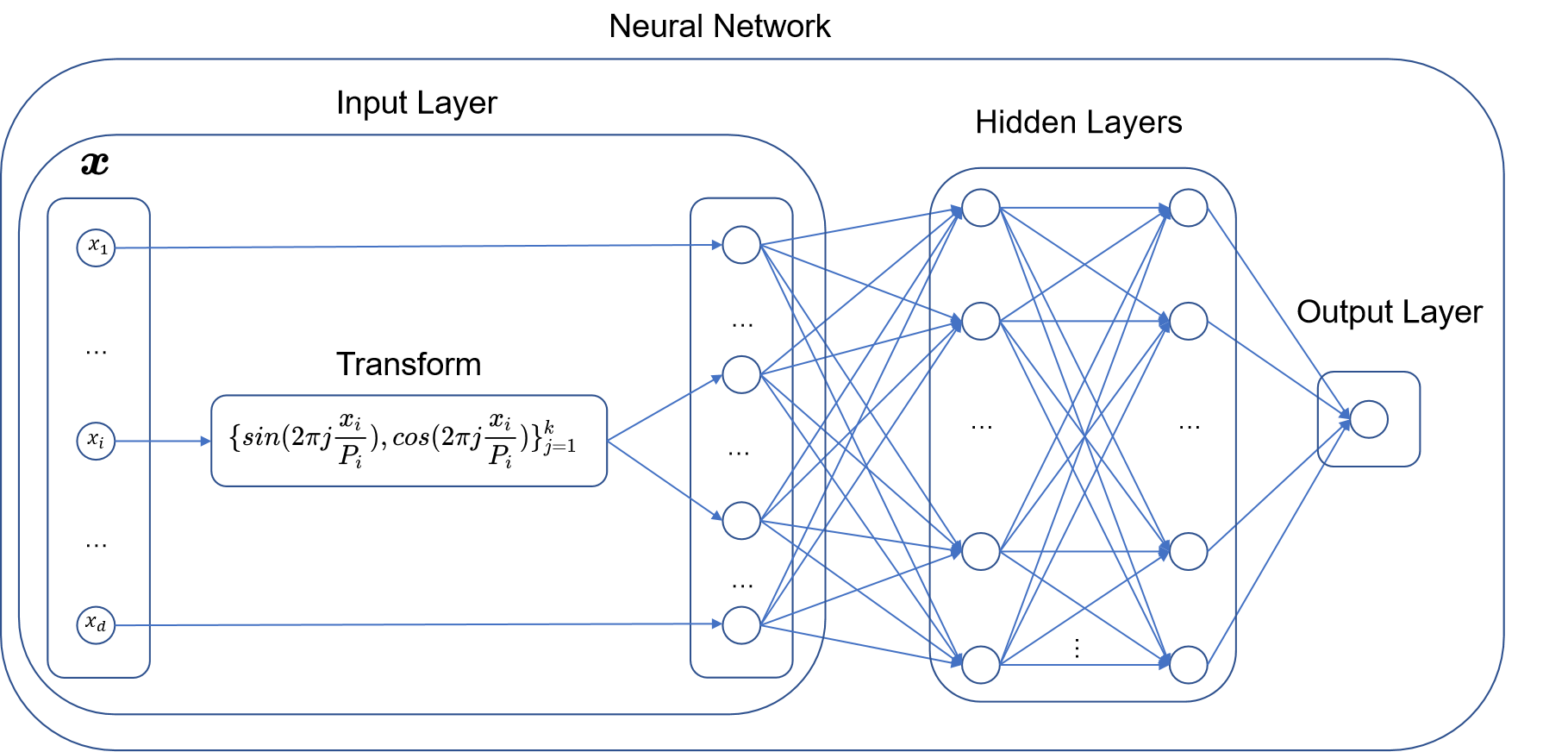}
    \caption{Network structure for periodic boundary condition.}
    \label{fig:u_enforce_periodic}
\end{figure}

We consider the periodic boundary condition of the form \eqref{eq:periodic_boundary_condition}, where $P_i$ is the period along the $i$-th direction. For more complex periodic boundary conditions, one can refer to \citep{dong2021method}.
\begin{equation}
    \label{eq:periodic_boundary_condition}
    \begin{array}{r@{}l}
        \begin{aligned}
             & u(x_1, \cdots, x_i+P_i, \cdots, x_d) = u(x_1, \cdots, x_i, \cdots, x_d), &  & \forall \boldsymbol{x} \in \partial \Omega.
        \end{aligned}
    \end{array}
\end{equation}
Unlike the enforcement of the Dirichlet boundary condition, which modify the output of neural networks, 
we need to modify the input before the first hidden layer of neural network. Suppose $u$ satisfies \eqref{eq:periodic_boundary_condition} in $x_i$ direction, then the component $x_i$ is transformed as follows,
\begin{equation}
    \label{eq:u_enforce_periodic}
    \begin{array}{r@{}l}
        \begin{aligned}
             & x_i \to \left\{ \sin\left(2\pi j \frac{x_i}{P_i}\right), \cos\left(2\pi j \frac{x_i}{P_i}\right) \right\}_{j=1}^k,
        \end{aligned}
    \end{array}
\end{equation}
where $k$ is the hyper-parameter. The network structure is shown in Figure \ref{fig:u_enforce_periodic}. 
 Obviously, the number of neurons of the neural network increases from $d$ to $(d+2k-1)$.

\section{Numerical Experiments}
\label{sec:experiments}
In this section, we present numerical results to demonstrate the applicability and accuracy of our methods, which include multi-dimensional simulations. The domain $\Omega$ with Dirichlet boundary condition is chosen as  $[0, 1]^d$. The domain $\Omega$ with period boundary condition is chosen as  $[0, 2\pi]^d$. In all the numerical experiments, we choose Adam optimizer with initial learning rate $10^{-3}$ to minimize the loss function, and we train the neural network with a simple architecture of MLP on a server equipped with CentOS 7 system, one Intel Xeon Platinum 8358 2.60GHz CPU and one NVIDIA A100 80GB GPU. Moreover, we choose $tanh$ as the active function and Latin hypercube sampling (lhs) \citep{loh1996latin} as the random sampling algorithm unless otherwise stated.
In addition, we train the neural network at a fixed number of epochs without selecting $\epsilon$ and the batch size equals to $N$. In the following numerical experiments, we use $\lVert u \rVert = \sqrt{\frac{1}{N}\sum_{i=1}^N u^2(\boldsymbol{x}_i)}$ to calculate the discretized norm and all the eigenfunctions that we show are normalized using $\lVert u \rVert$. It is hard to visualize the high-dimensional eigenfunction directly. Therefore, the density of a function $u$ is defined as the probability density function of $u(\mathbf{X})$, where $\mathbf{X}$ is a uniformly distributed random variable on $\Omega$.


\subsection{An Example of Eigenvalue Problem implemented by PMNN}
To validate the performance of PMNN, we solve the following problem in $\Omega=[0, 1]^d$ and $d\leq 10$, 
\begin{equation}
    \label{eq:dominant_problem}
    \begin{array}{r@{}l}
        \left\{
        \begin{aligned}
            \Delta u + 100 u & =\lambda u, &  & \mbox{in} \enspace \Omega,          \\
            u             & = 0,        &  & \mbox{on} \enspace \partial \Omega.
        \end{aligned}
        \right.
    \end{array}
\end{equation}
The exact solution of the largest positive eigenvalue is given by $\lambda=100-d\pi^2$ for $d\leq 10$ and the associated eigenfunction is $u=\Pi_{i=1}^d \sin(\pi x_i)$. The parameters used to train PMNN for different dimensions are summarized in Table \ref{tab:params_pmnn}. We choose the number of random sampling points $N=10000, 20000, 50000, 100000$ for $d=1, 2, 5, 10$, respectively. In high dimensions, we train PMNN with more epochs and more neurons.

\begin{figure}[htbp]
    \begin{minipage}{0.48\linewidth}
        \centering
        \includegraphics[width=1\textwidth]{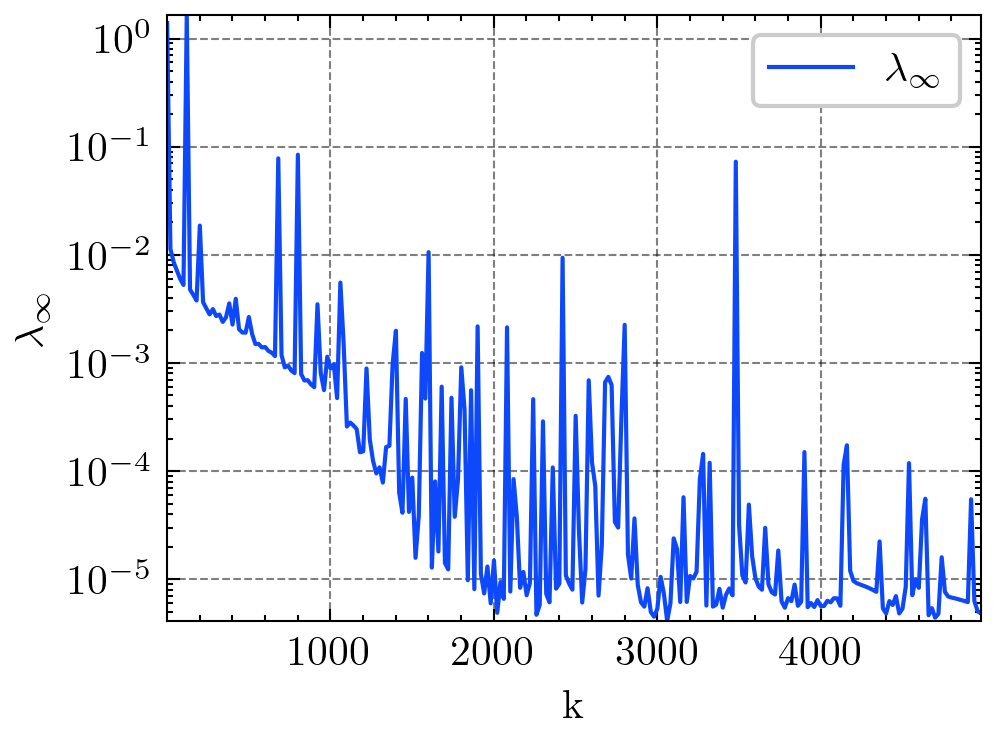}
    \end{minipage}
    \begin{minipage}{0.48\linewidth}
        \centering
        \includegraphics[width=1\textwidth]{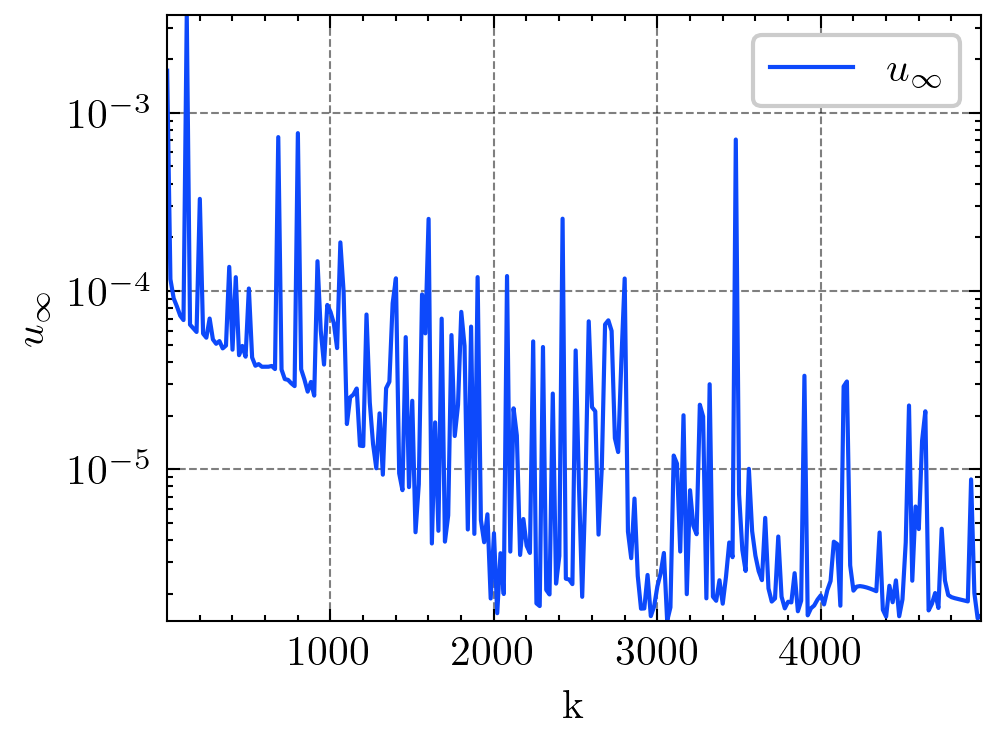}
    \end{minipage}
    \caption{The maximum norm of the eigenvalue and the associated eigenfunction of the problem \eqref{eq:dominant_problem} with iteration $k$ increasing in $d=1$.}
    \label{fig:pmnn_convergence}
\end{figure}

\begin{figure}[htbp]
    \begin{minipage}{0.48\linewidth}
        \centering
        \includegraphics[width=1\textwidth]{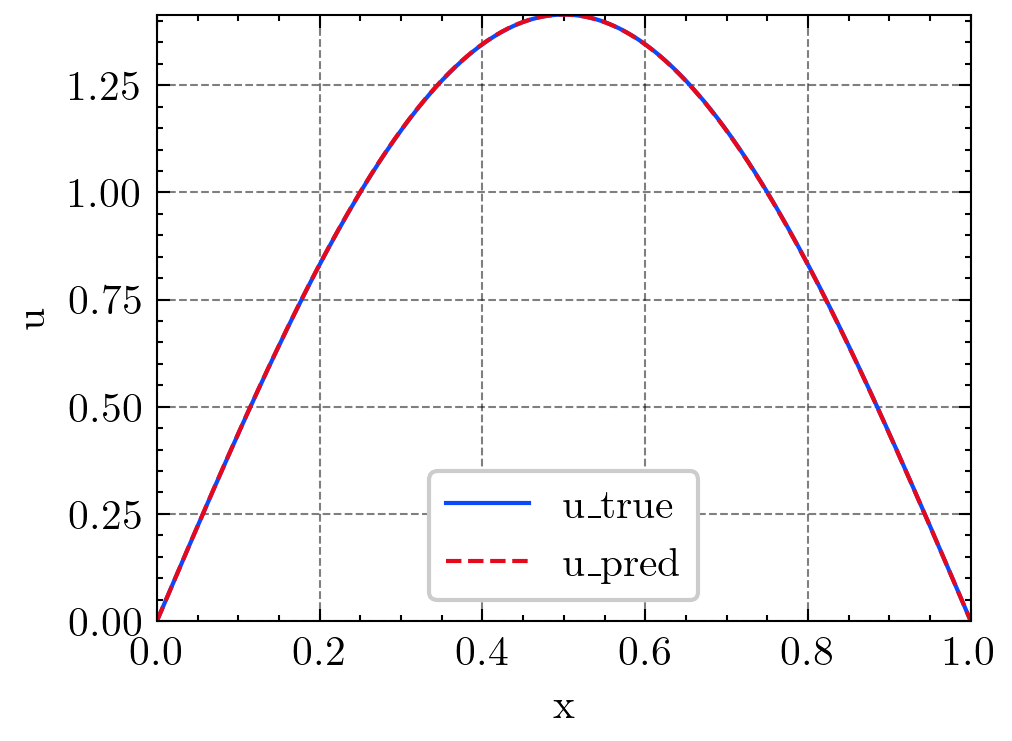}
    \end{minipage}
    \begin{minipage}{0.48\linewidth}
        \centering
        \includegraphics[width=1\textwidth]{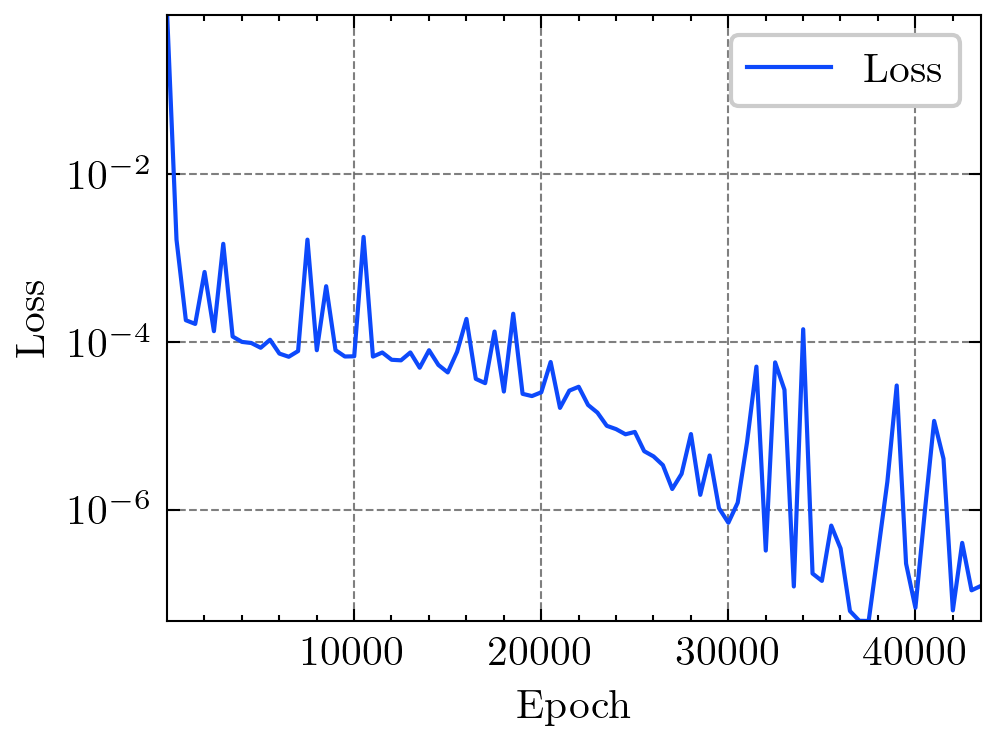}
    \end{minipage}
    \caption{Left: the eigenfunction of \eqref{eq:dominant_problem} learned by PMNN  and the exact eigenfunction  in one dimension. Right: loss in the training process.}
    \label{fig:pmnn_1D_line}
\end{figure}

\begin{figure}[htp]
    \centering
    \includegraphics[width=1\textwidth]{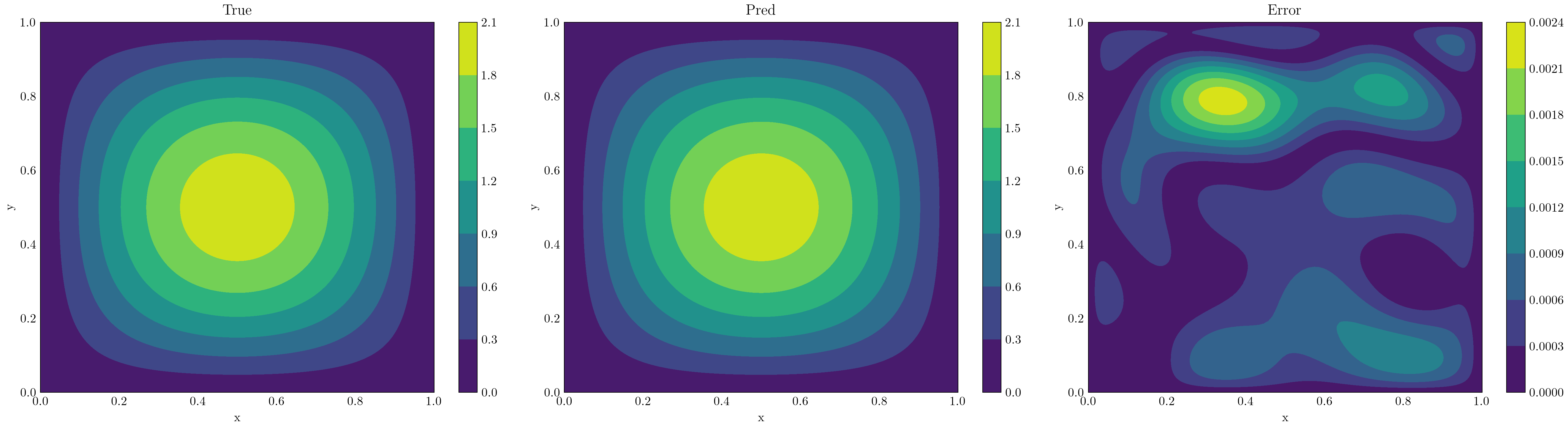}
    \caption{Heat maps of the eigenfunction of \eqref{eq:dominant_problem} in two dimension. Left: the exact solution; Middle: the prediction solution computed by PMNN; Right: the absolute error.}
    \label{fig:pmnn_2D_heatmap3}
\end{figure}

\begin{figure}[htbp]
    \begin{minipage}{0.48\linewidth}
        \centering
        \includegraphics[width=1\textwidth]{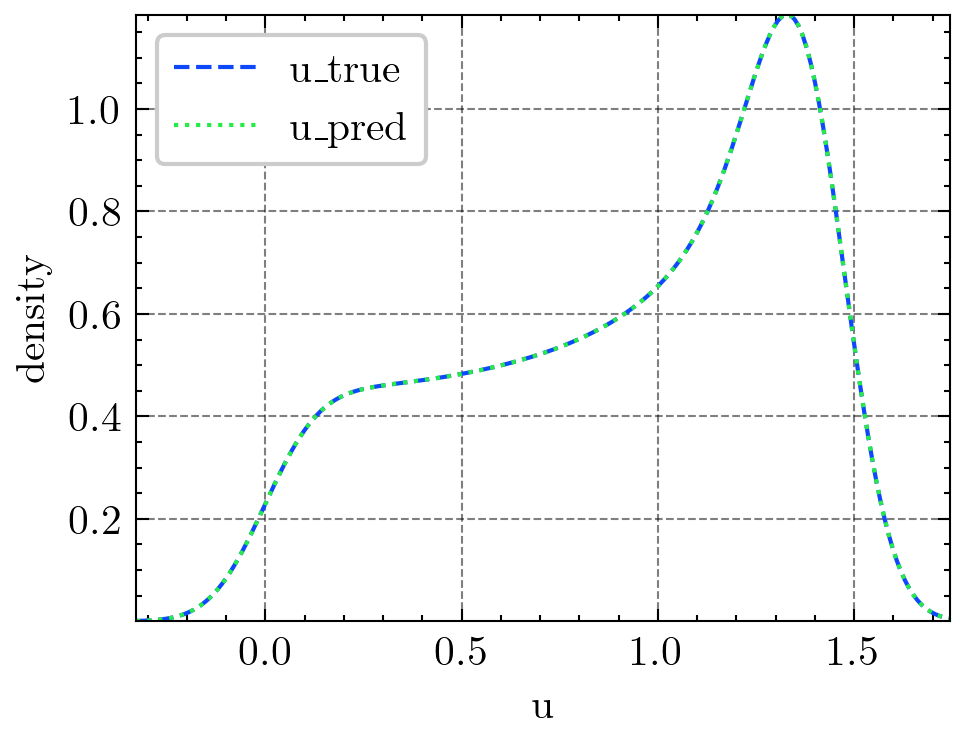}
        \label{fig:pm_density_1}
    \end{minipage}
    \begin{minipage}{0.48\linewidth}
        \centering
        \includegraphics[width=1\textwidth]{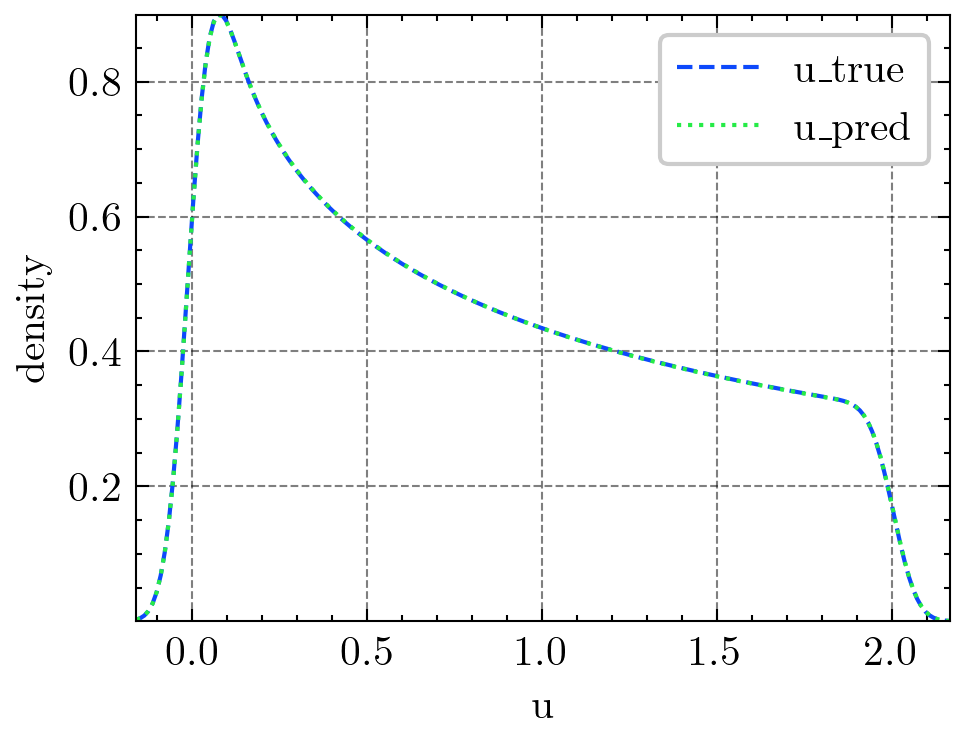}
        \label{fig:pm_density_2}
    \end{minipage}
    \begin{minipage}{0.48\linewidth}
        \centering
        \includegraphics[width=1\textwidth]{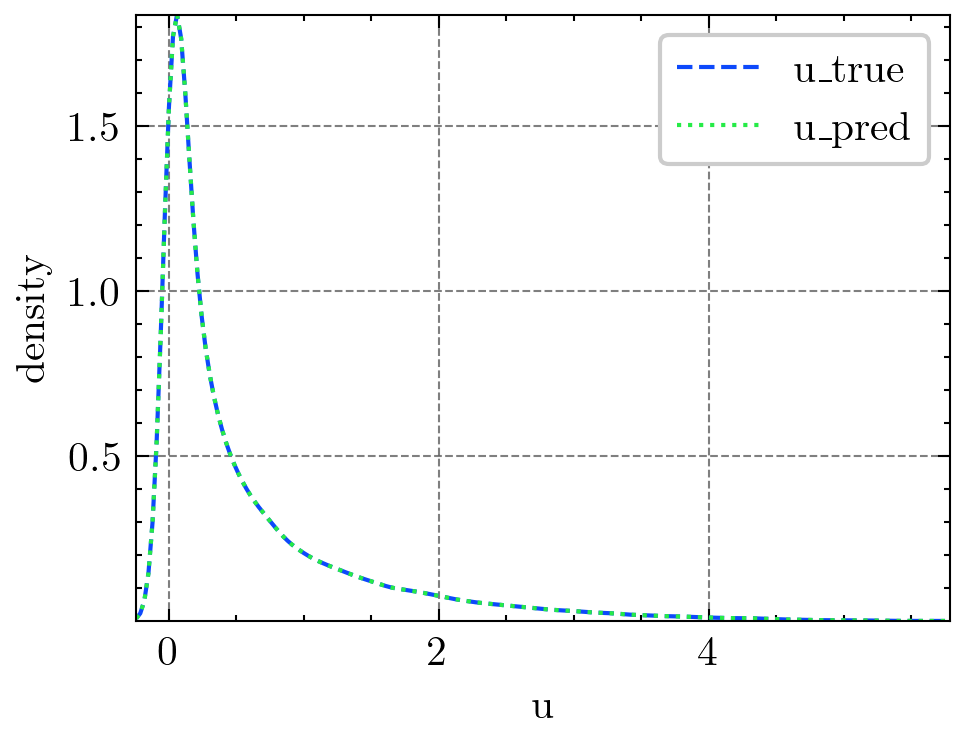}
        \label{fig:pm_density_5}
    \end{minipage}
    \begin{minipage}{0.48\linewidth}
        \centering
        \includegraphics[width=1\textwidth]{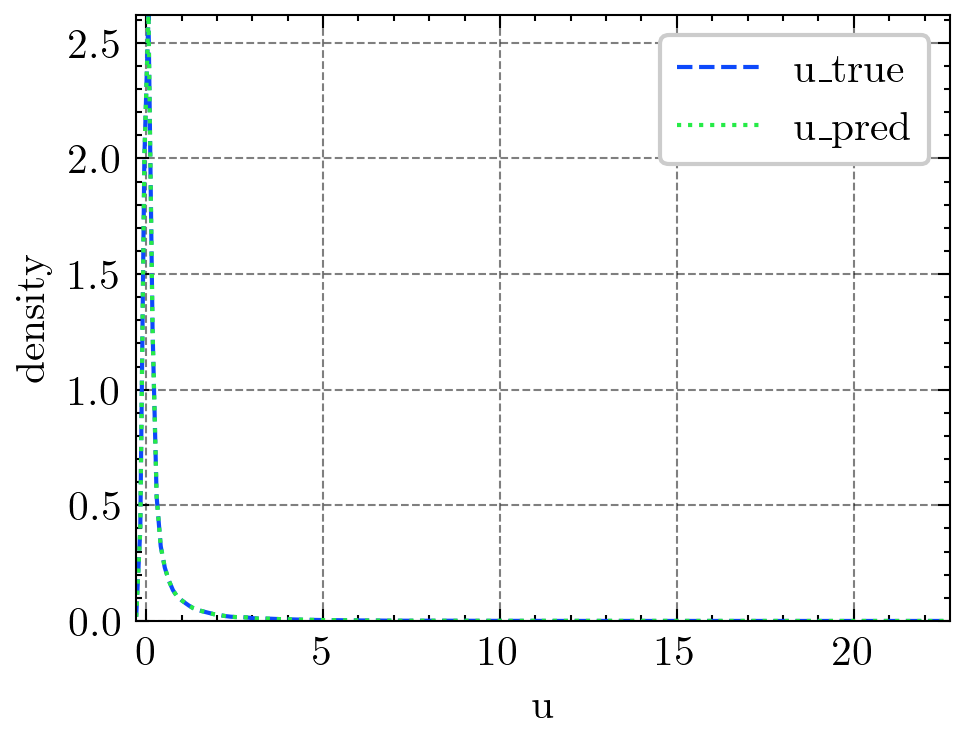}
        \label{fig:pm_density_10}
    \end{minipage}
    \caption{For equation \eqref{eq:dominant_problem}, comparison of the densities of approximate eigenfunctions learned by PMNN  and the exact eigenfunctions in (Top:) $d = 1$ (left) and $d= 2$ (right); (Bottom:) $d=5$ (left)  and $d=10$ (right).}
    \label{fig:pm_density}
\end{figure}

We use $\lambda_{\infty}$ to denote $\max (\lvert \lambda_{pred} - \lambda_{true} \lvert)$ and $u_{\infty}$ to denote $L_\infty$ norm of u. The maximum norm of $\lambda_{\infty}$ and $u_{\infty}$ of the eigenvalue problem \eqref{eq:dominant_problem} with iteration $k$ increasing (training process proceeding) in $d=1$ is shown in Figure \ref{fig:pmnn_convergence}.  It is observed that the maximum norm goes to zero with some variations, which implies the convergence of PMNN.
In Figure \ref{fig:pmnn_1D_line}, the eigenfunction obtained by PMNN, the exact eigenfunction and the variation of loss are presented. It is easy to see that PMNN perfectly learns the eigenfunction in one dimensional case and the loss decreases in the training process. 
The eigenfunction learned by PMNN in two dimension, the exact solution and the absolute error between the NN solution and the true solution are shown in Figure \ref{fig:pmnn_2D_heatmap3}. It is obvious to see that PMNN perfectly learns the eigenfunction in two dimension.
To further validate the performance of our method, we also implement our method in higher dimensions and the results are shown in Table \ref{tab:error_pmnn} and Figure \ref{fig:pm_density}. From Table \ref{tab:error_pmnn}, the relative error is small enough to demonstrate the accuracy of PMNN in different dimensions. Figure \ref{fig:pm_density} shows the densities of eigenfunctions of equation \eqref{eq:dominant_problem} in $d= 1$, $d = 2$, $d = 5$ and $d =10$. It is shown that the densities of eigenfunctions learned by PMNN perfectly fit the densities of the exact eigenfunctions in all cases.

\begin{table}[htp]
\begin{center}
    \caption{Parameter settings of training PMNN in different dimensions.}
    \begin{tabular}{llll}
        \hline\noalign{\smallskip}
          $d$ & $N$    & $N_{epoch}$ & layers of MLP\\
        \hline
         1   & 10000  & 50000       & [1, 20, 20, 20, 20, 1]  \\
         2   & 20000  & 50000       & [2, 20, 20, 20, 20, 1]  \\
         5   & 50000  & 50000       & [5, 40, 40, 40, 40, 1]  \\
         10  & 100000 & 100000      & [10, 80, 80, 80, 80, 1] \\
        \hline
    \end{tabular}
    \label{tab:params_pmnn}
    \end{center}
    \end{table}

\begin{table}[htp]
\begin{center}
    \caption{Comparison of exact eigenvalues and approximate eigenvalues learned by PMNN in different dimensions.}
    \begin{tabular}{lllll}
        \hline\noalign{\smallskip}
          $d$ & Exact $\lambda$ & Approximate $\lambda$  &  Relative error \\
        \hline
         1   & 90.1304         & 90.1302                   & 2.4727E-06     \\
         2   & 80.2608         & 80.2603                    & 6.1240E-06     \\
         5   & 50.6520         & 50.6513                    & 1.3264E-05     \\
         10  & 1.3040          & 1.3004                     & 2.7461E-03     \\      
        \hline
    \end{tabular}
    \label{tab:error_pmnn}
    \end{center}
\end{table}


\subsection{Examples of Smallest Eigenvalue Problem}
\subsubsection{Harmonic Eigenvalue Problem}
\begin{figure}[htbp]
    \begin{minipage}{0.48\linewidth}
        \centering
        \includegraphics[width=1\textwidth]{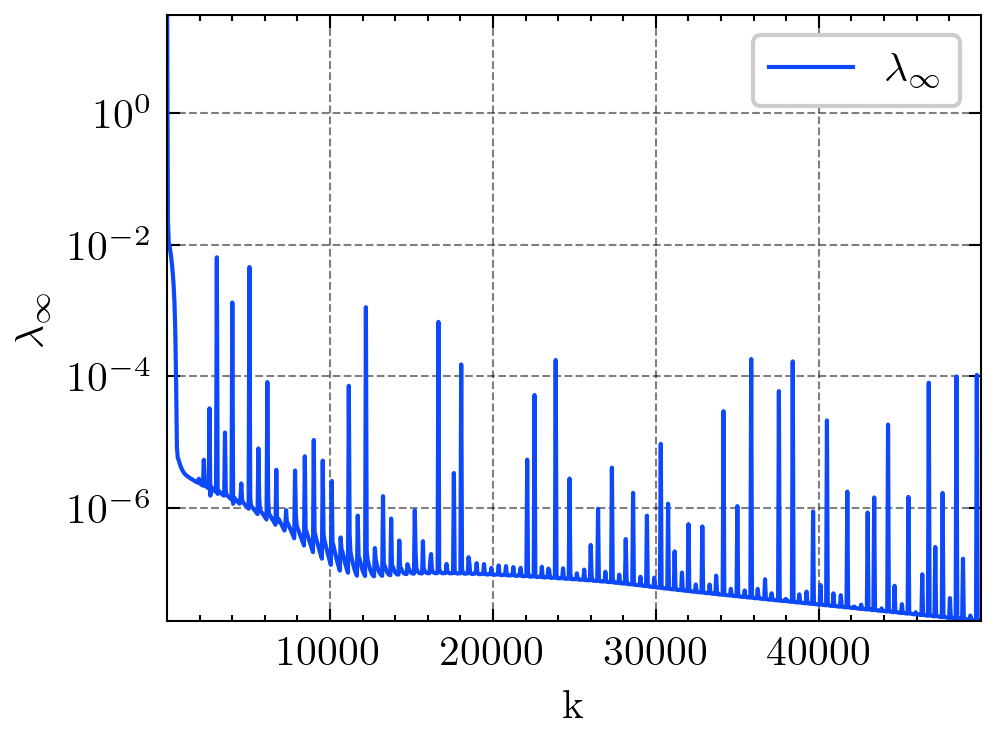}
    \end{minipage}
    \begin{minipage}{0.48\linewidth}
        \centering
        \includegraphics[width=1\textwidth]{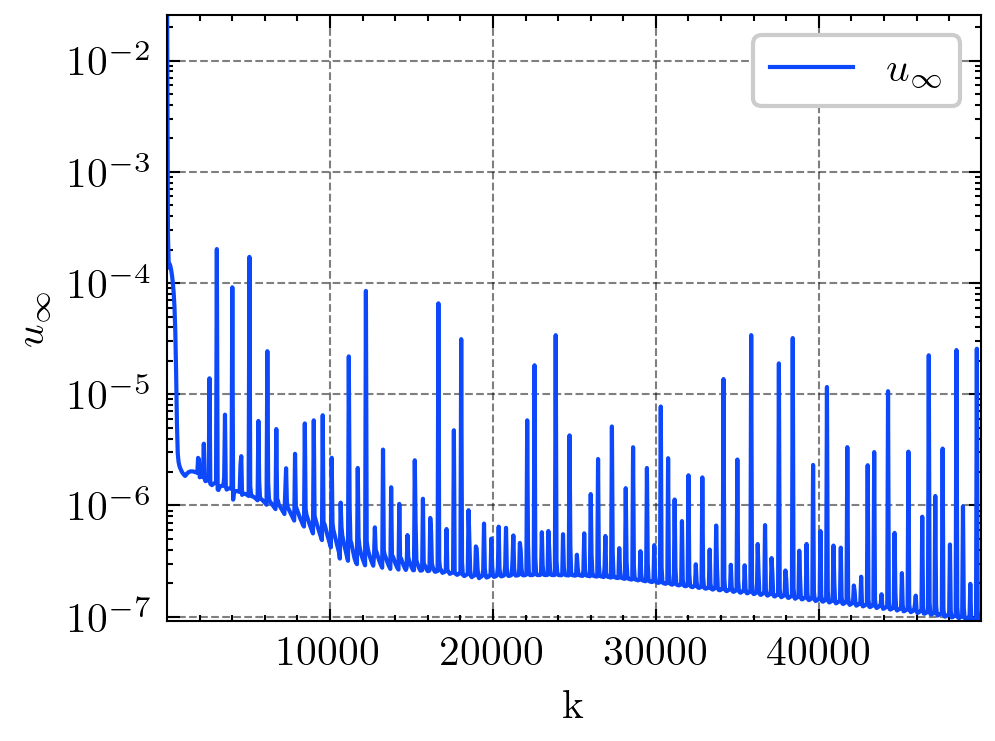}
    \end{minipage}
    \caption{The maximum norm of the eigenvalue and the associated eigenfunction of the problem \eqref{eq:smallest_problem1} with iteration $k$ increasing in $d=1$.}
    \label{fig:ipmnn_convergence}
\end{figure}

\begin{figure}[htbp]
     \begin{minipage}{0.48\linewidth}
        \centering
        \includegraphics[width=1\textwidth]{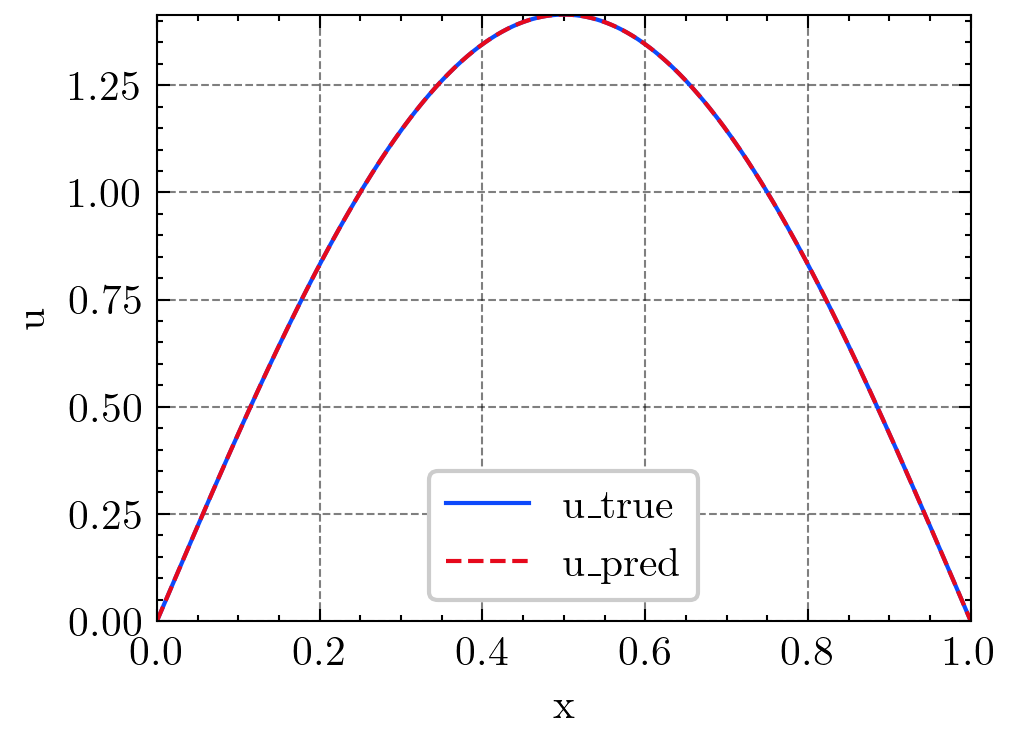}
    \end{minipage}
    \begin{minipage}{0.48\linewidth}
        \centering
        \includegraphics[width=1\textwidth]{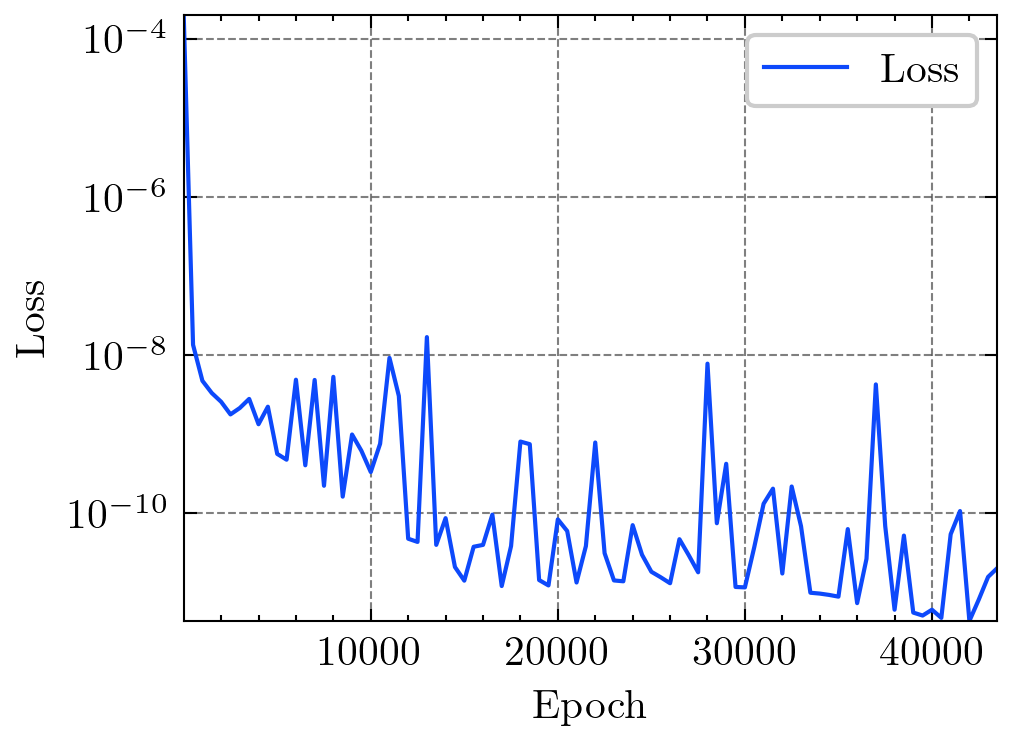}
        \label{fig:ipmnn_1D_loss}
    \end{minipage}
    \caption{Left: the eigenfunction of \eqref{eq:smallest_problem1} learned by IPMNN  and the exact eigenfunction in one dimension. Right: loss in the training process.}    \label{fig:ipmnn_1D_line}
\end{figure}

\begin{figure}[htp]
    \centering
    \includegraphics[width=1\textwidth]{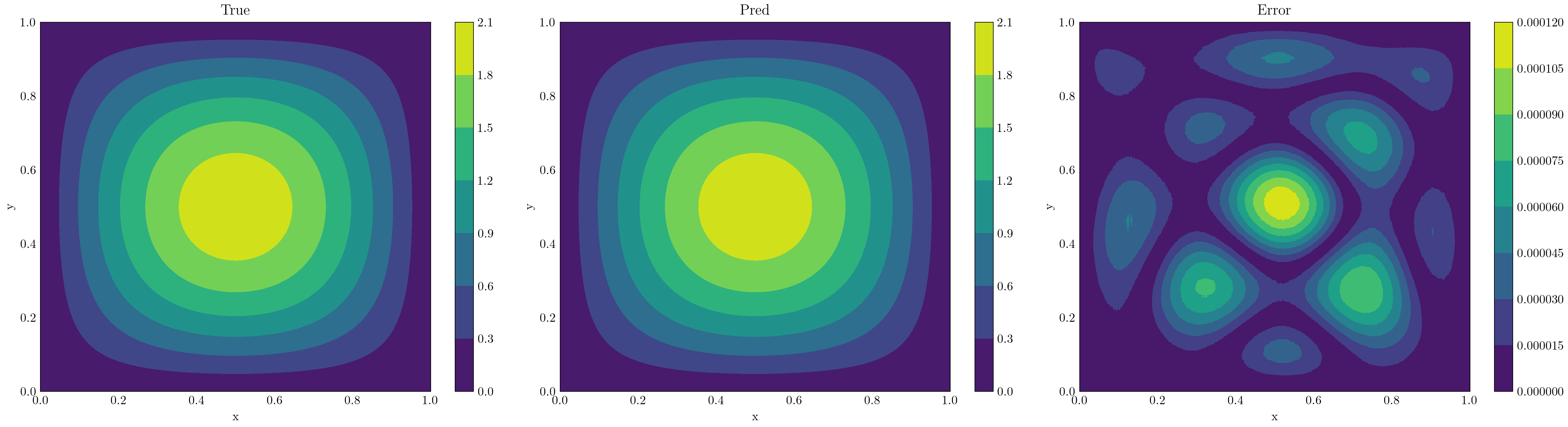}
    \caption{Heat maps of the eigenfunction of \eqref{eq:smallest_problem1} in two dimension. Left: the exact solution; Middle: the prediction solution computed by IPMNN; Right: the absolute error.}
    \label{fig:ipmnn_2D_heatmap3}
\end{figure}

\begin{figure}[htbp]
    \begin{minipage}{0.48\linewidth}
        \centering
        \includegraphics[width=1\textwidth]{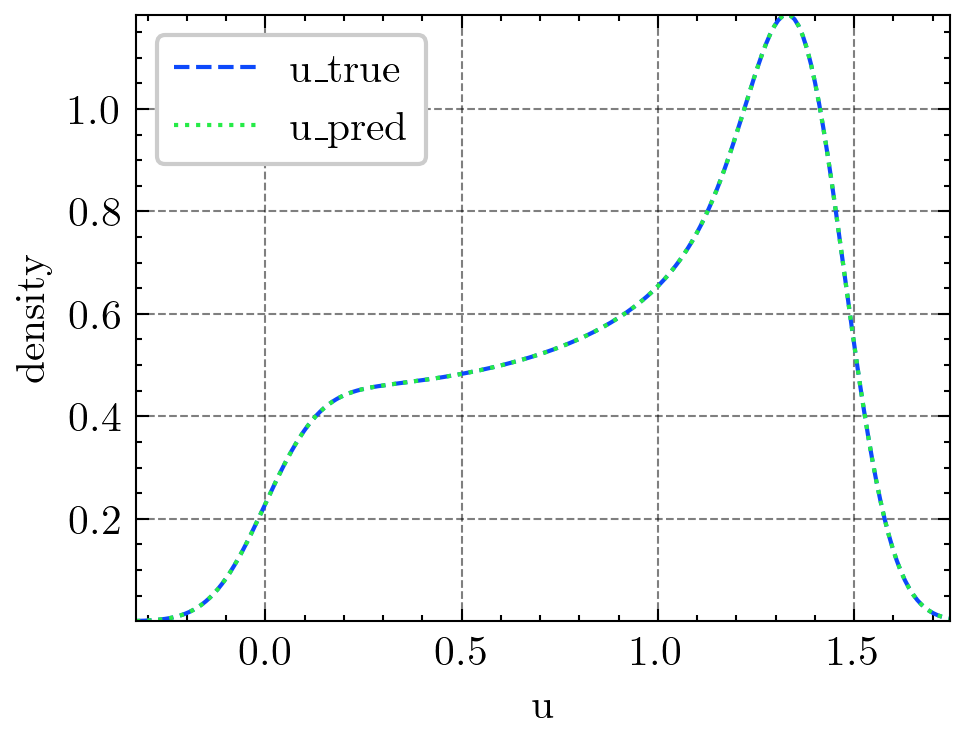}
        \label{fig:ipm_density_1}
    \end{minipage}
    \begin{minipage}{0.48\linewidth}
        \centering
        \includegraphics[width=1\textwidth]{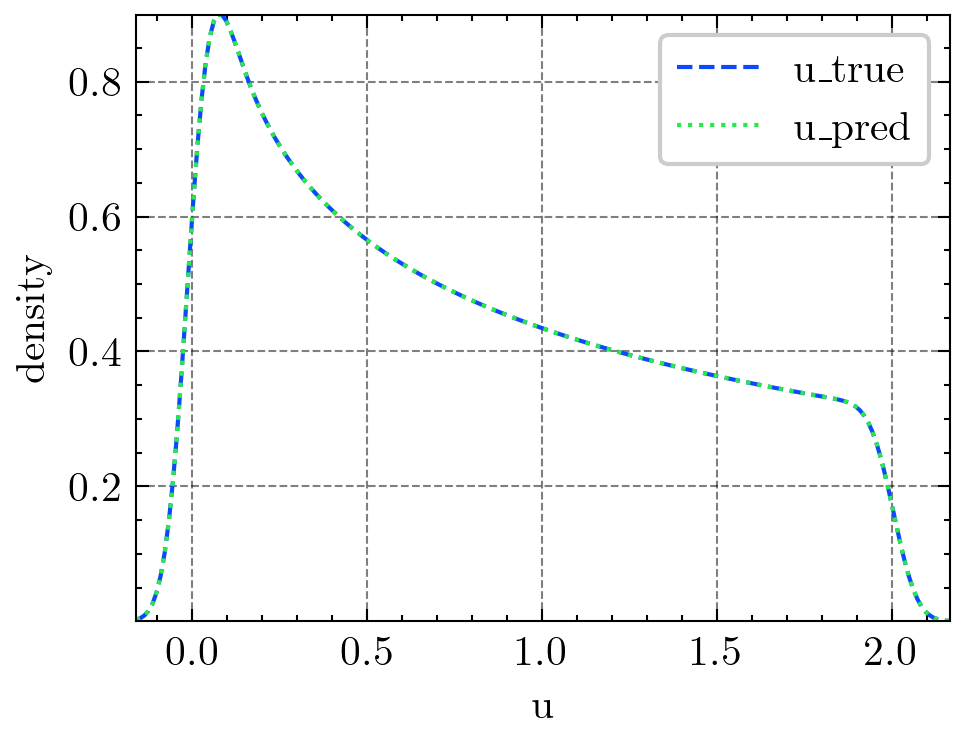}
        \label{fig:ipm_density_2}
    \end{minipage}
    \begin{minipage}{0.48\linewidth}
        \centering
        \includegraphics[width=1\textwidth]{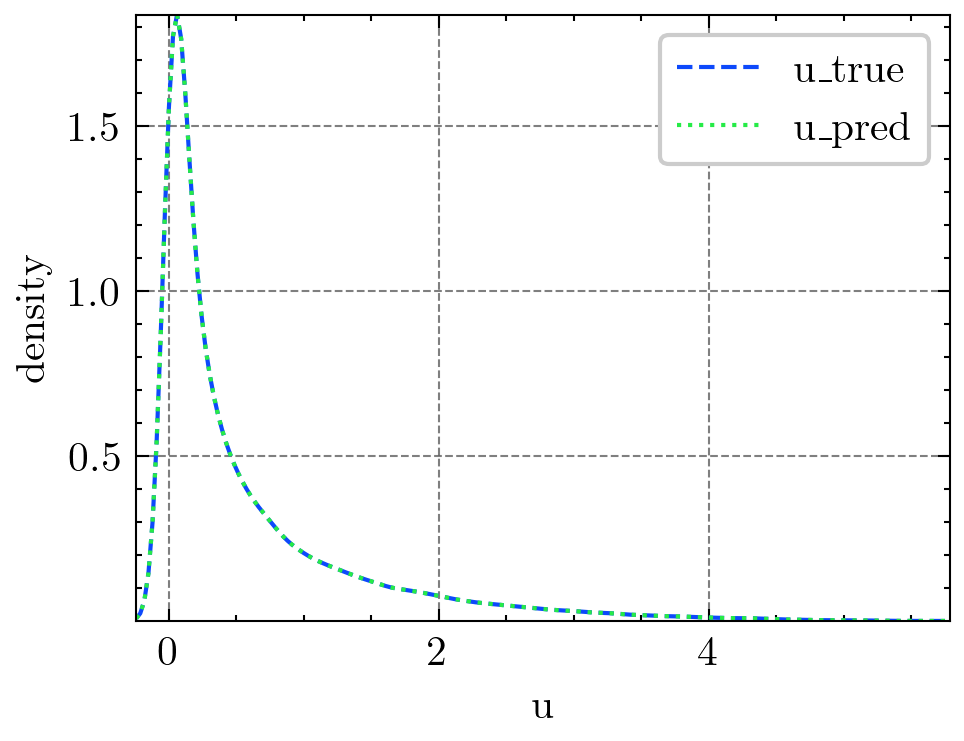}
        \label{fig:ipm_density_5}
    \end{minipage}
    \begin{minipage}{0.48\linewidth}
        \centering
        \includegraphics[width=1\textwidth]{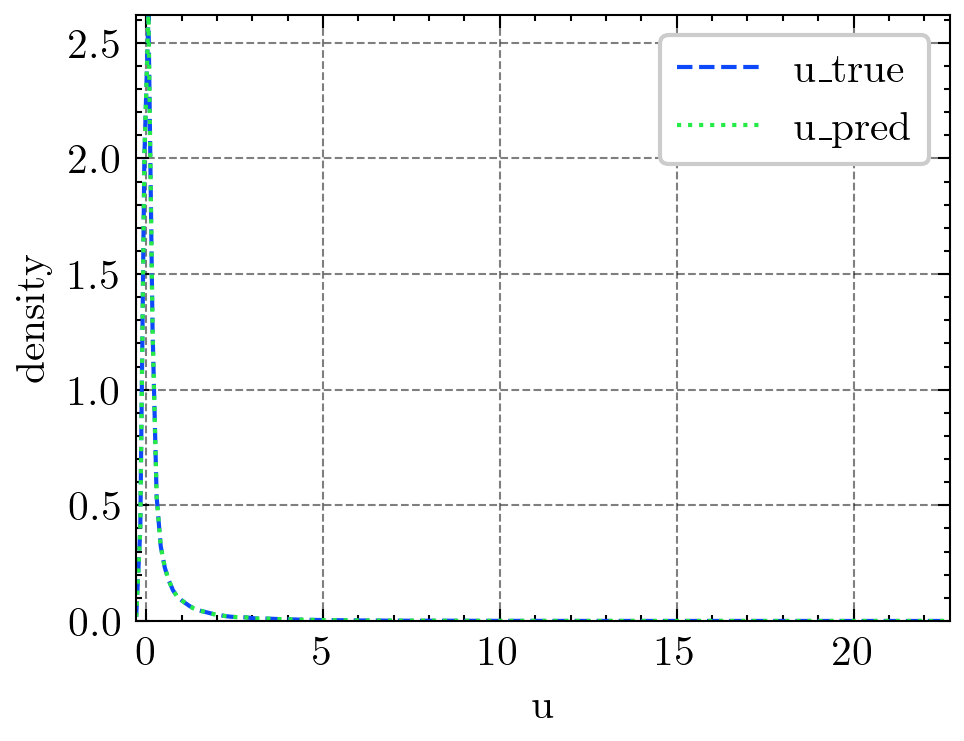}
        \label{fig:ipm_density_10}
    \end{minipage}
    \caption{For equation \eqref{eq:smallest_problem1}, comparison of the densities of approximate eigenfunctions learned by IPMNN  and  the exact eigenfunctions in (Top:) $d = 1$ (left) and $d= 2$ (right); (Bottom:) $d=5$ (left)  and $d=10$ (right). }   
    \label{fig:ipm_density}
\end{figure}

To validate the performance of IPMNN, we solve the following problem in $\Omega=[0, 1]^d$,
\begin{equation}
    \label{eq:smallest_problem1}
    \begin{array}{r@{}l}
        \left\{
        \begin{aligned}
            -\Delta u & =\lambda u, &  & \mbox{in} \enspace \Omega,          \\
            u         & = 0,        &  & \mbox{on} \enspace \partial \Omega.
        \end{aligned}
        \right.
    \end{array}
\end{equation}
The smallest eigenvalue is given by $\lambda=d\pi^2$ and the associated eigenfunction is $u=\Pi_{i=1}^d\sin(\pi x_i)$.  We use the same parameters as in Table \ref{tab:params_pmnn}. 
The maximum norm of $\lambda_{\infty}$ and $u_{\infty}$ of the eigenvalue problem \eqref{eq:smallest_problem1} with iteration $k$ increasing (training process proceeding) in $d=1$ is shown in Figure \ref{fig:ipmnn_convergence}. Obviously, IPMNN is able to find the smallest eigenvalue and the associated eigenfunction.
In Figure \ref{fig:ipmnn_1D_line}, the smallest eigenfunction obtained by IPMNN, the exact eigenfunction  and the variation of loss are presented. It is easy to see that IPMNN perfectly learns the eigenfunction in one dimensional case and the loss decreases in the training process. The eigenfunction learned by IPMNN in two dimension, the exact solution and the absolute error between the NN solution and the true solution are shown in Figure \ref{fig:ipmnn_2D_heatmap3}. It is obvious to see that IPMNN perfectly learns the eigenfunction in two dimension. We also implement our method in higher dimensions and we compare our results with that obtained by DRM \citep{yu2018deep} in Table \ref{tab:error_ipmnn}. The relative error is small enough to demonstrate the accuracy of IPMNN in different dimensions. The densities of eigenfunctions of equation\eqref{eq:smallest_problem1} in $d= 1$, $d = 2$, $d = 5$ and $d =10$ are shown in  Figure \ref{fig:ipm_density}. It is obvious to see that the densities of eigenfunctions learned by IPMNN perfectly fit the densities of the exact eigenfunctions in all cases.



\begin{table}[htp]
\begin{center}
    \caption{Comparison of exact eigenvalues and approximate eigenvalues for \eqref{eq:smallest_problem1} in different dimensions.}
    \begin{tabular}{lllll}
        \hline\noalign{\smallskip}
        %
           Method                      & $d$ & Exact $\lambda$ & Approximate $\lambda$  &  Relative error \\
        \hline
        IPMNN                        & 1   & 9.8696          & 9.8695                     & 7.6630E-06     \\
        Deep Ritz \citep{yu2018deep} & 1   & 9.8696          & 9.85                                & 2.0E-03        \\
        IPMNN                        & 2   & 19.7392         & 19.7395                    & 1.3209E-05     \\
        IPMNN                        & 5   & 49.3480         & 49.3485                    & 9.8266E-06     \\
        Deep Ritz \citep{yu2018deep} & 5   & 49.3480         & 49.29                               & 1.1E-03        \\
        IPMNN                        & 10  & 98.6960         & 98.6953                    & 7.9529E-06     \\
        Deep Ritz \citep{yu2018deep} & 10  & 98.6960         & 92.35                               & 6.43E-02       \\
        \hline
    \end{tabular}
    \label{tab:error_ipmnn}
\end{center}
\end{table}


\subsubsection{Comparison with Finite Difference Method}

We compare the numerical results obtained by IPMNN and FDM. In $d =2$, we use the inverse power method to solve for the smallest eigenvalue problem \eqref{eq:smallest_problem1}, which discretized in a uniform grid and computed by FDM. For comparison, the same uniform sampling points are utilized to train IPMNN.

\begin{figure}[htbp]
     \begin{minipage}{0.9\linewidth}
        \centering
        \includegraphics[width=1\textwidth]{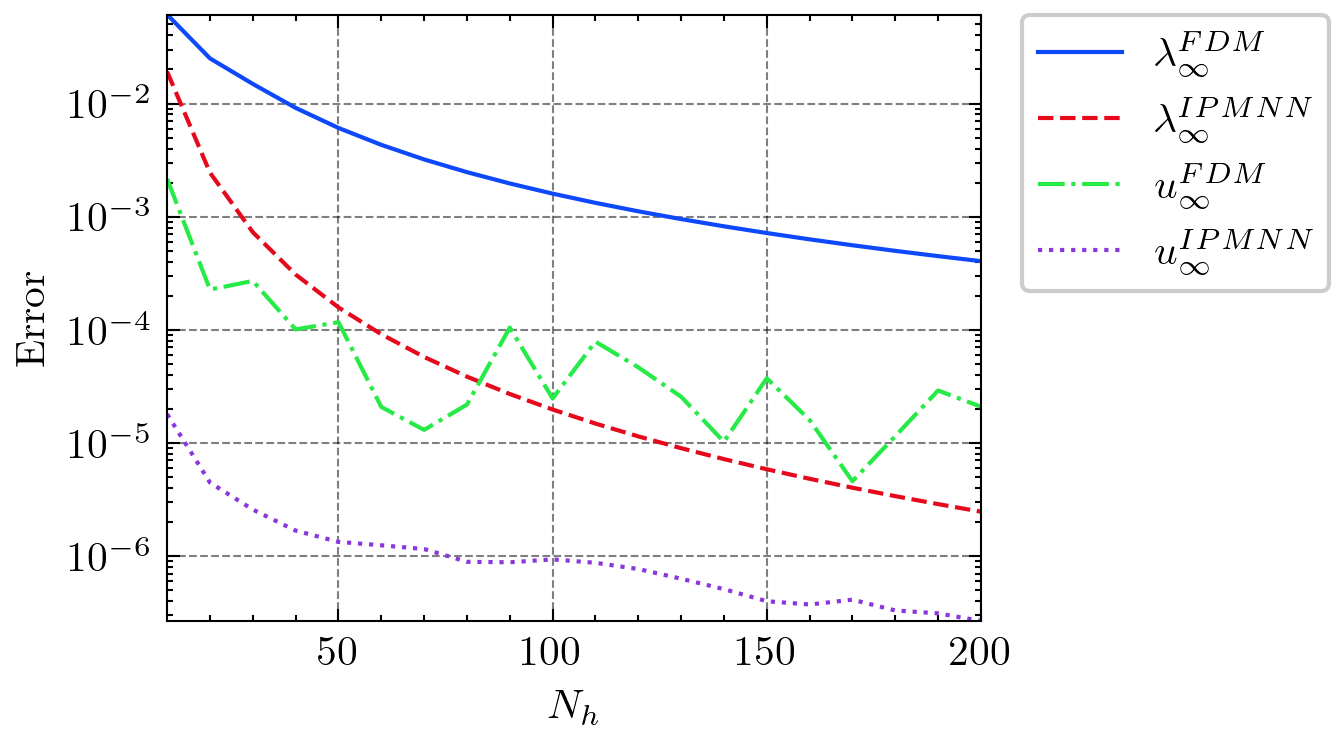}
    \end{minipage}
    \caption{Comparison of the errors numerically obtained by FDM and PINN in 2D.}
    \label{fig:FDM_IPMNN_abs}
\end{figure}

In Figure \ref{fig:FDM_IPMNN_abs}, the horizontal axis denotes the number of points $N_h$ on the x or y axis in computational domain and the total number of points is $N_h^2$. The vertical axis represents the maximum norm of the errors of the eigenvalue and the associated eigenfunction. For both FDM and IPMNN, the accuracy of the eigenvalue and the associated eigenfunction increases with the number of training points. It is observed, IPMNN is better than FDM in its ability to approximate the eigenvalue and the associated eigenfunction. However, it should be noted that much more training time for the neural network is needed than the cost for FDM. In addition, the traditional numerical methods are difficult to solve high-dimensional problems and IPMNN is suitable for high-dimensional problems as shown in our experiments.

\subsubsection{Fokker-Planck Equation}
\begin{figure}[htbp]
    \begin{minipage}{0.48\linewidth}
        \centering
        \includegraphics[width=1\textwidth]{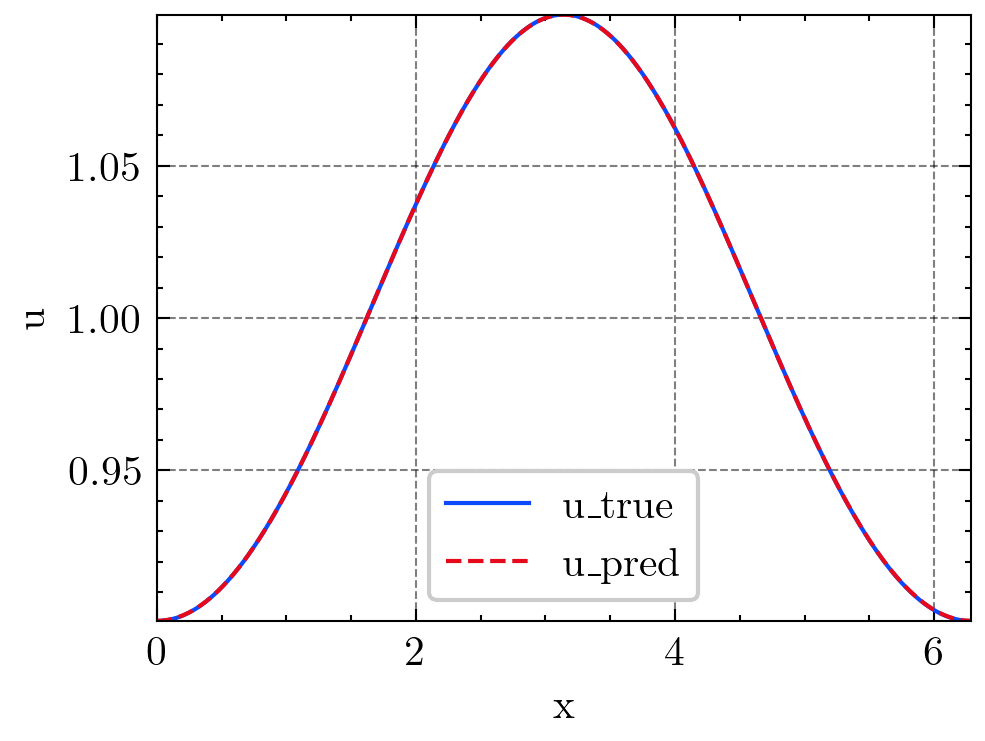}
    \end{minipage}
    \begin{minipage}{0.48\linewidth}
        \centering
        \includegraphics[width=1\textwidth]{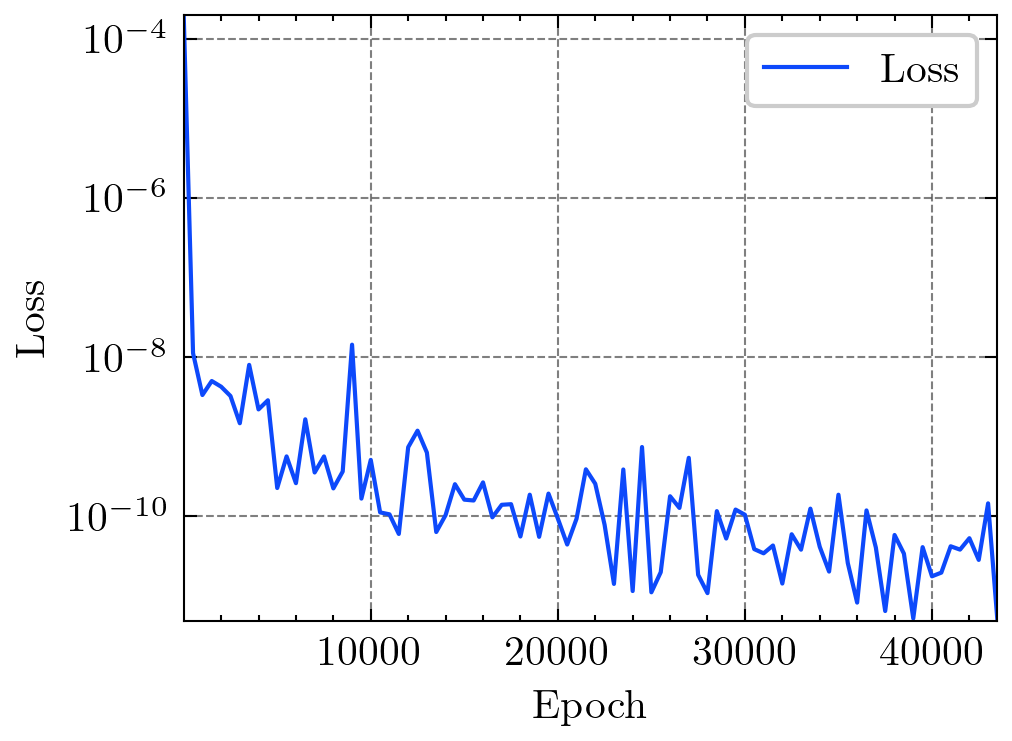}
        \label{fig:f_loss}
    \end{minipage}
     \caption{Left: the eigenfunction of Fokker-Planck equation \eqref{eq:fokker_planck_equation} learned by IPMNN and the exact eigenfunction in one dimension. Right: loss in the training process.}
    \label{fig:fk_line}
\end{figure}

\begin{figure}[htp]
    \centering
    \includegraphics[width=1\textwidth]{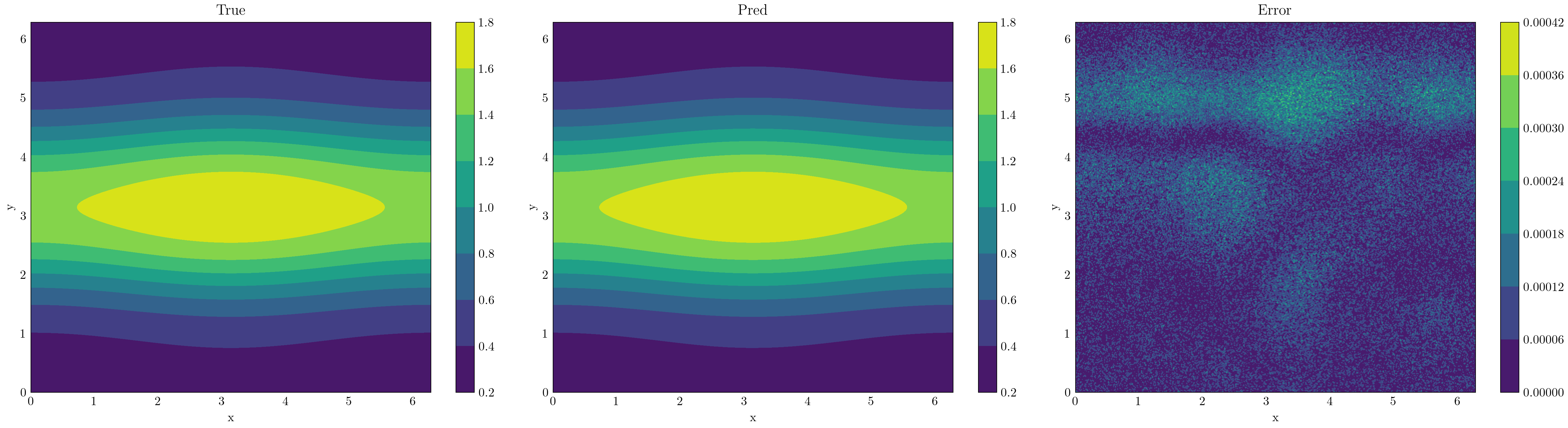}
    \caption{Heat maps of the eigenfunction of  Fokker-Planck equation \eqref{eq:fokker_planck_equation} in two dimension. Left: the exact solution; Middle: the prediction solution computed by IPMNN; Right: the absolute error.}
    \label{fig:ipmnn_2D_heatmap3_fk}
\end{figure}

\begin{figure}[htbp]
    \begin{minipage}{0.48\linewidth}
        \centering
        \includegraphics[width=1\textwidth]{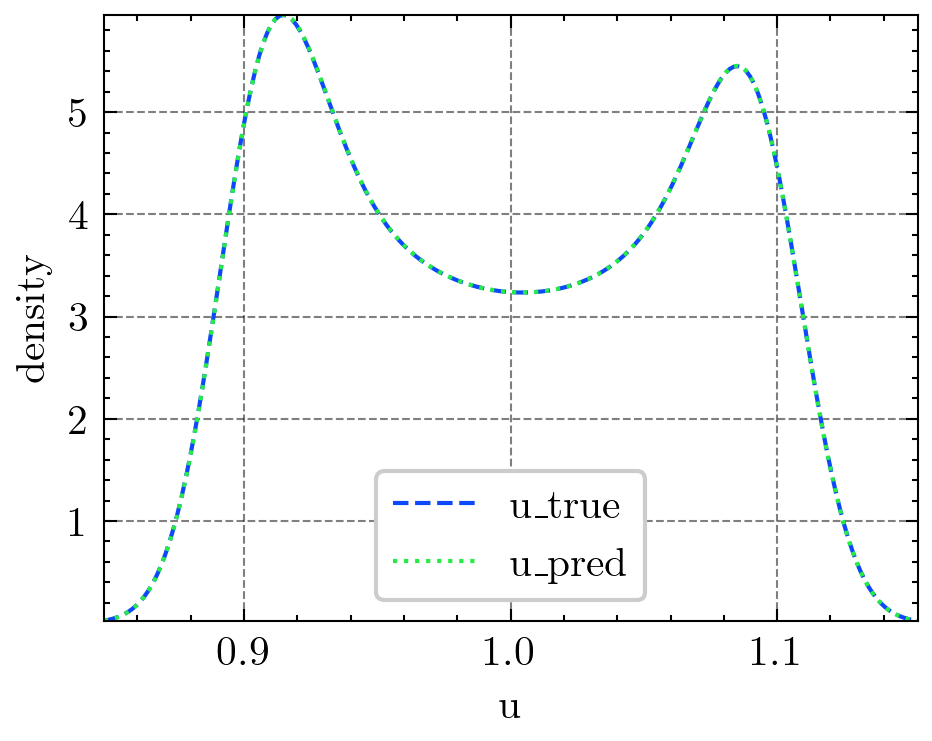}
        \label{fig:fk_density_1}
    \end{minipage}
    \begin{minipage}{0.48\linewidth}
        \centering
        \includegraphics[width=1\textwidth]{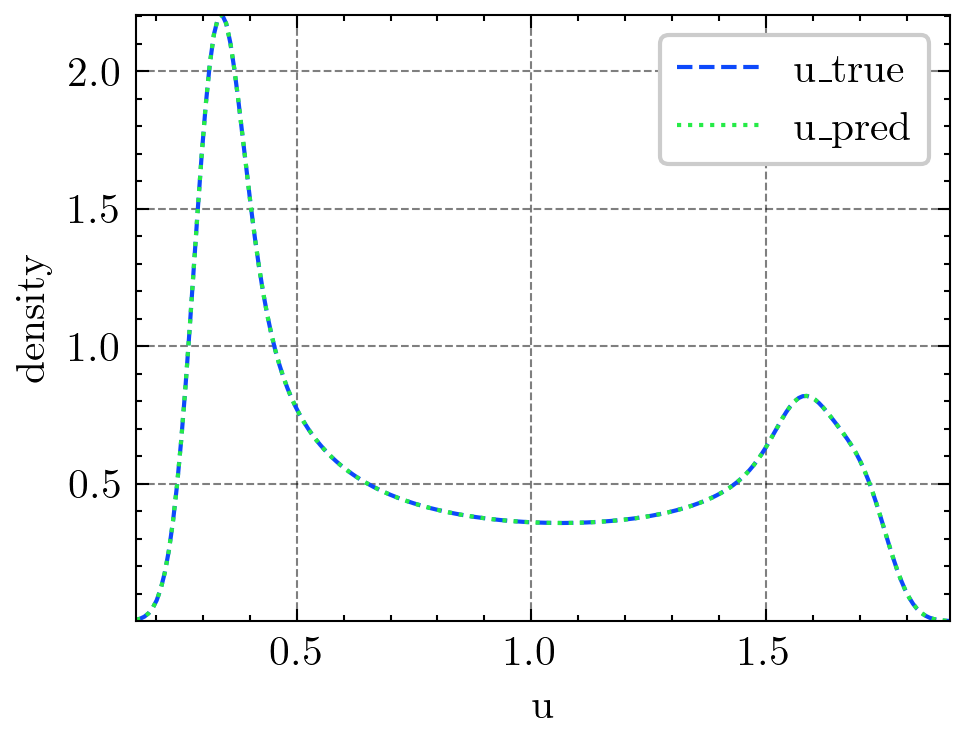}
        \label{fig:fk_density_2}
    \end{minipage}
    \begin{minipage}{0.48\linewidth}
        \centering
        \includegraphics[width=1\textwidth]{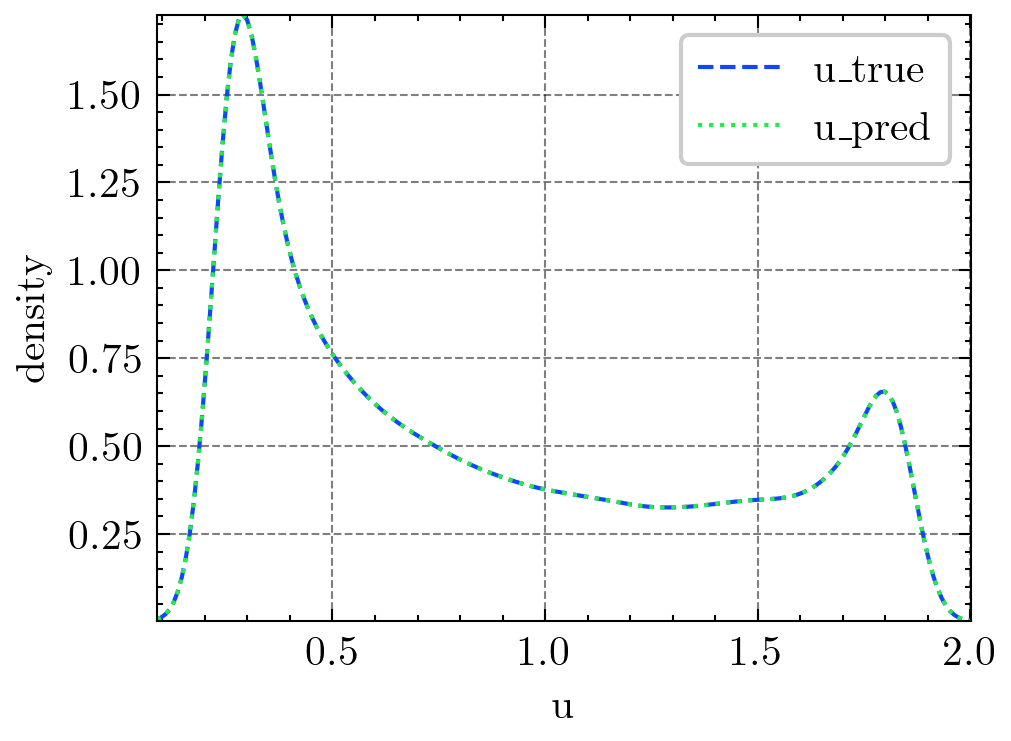}
        \label{fig:fk_density_5}
    \end{minipage}
    \begin{minipage}{0.48\linewidth}
        \centering
        \includegraphics[width=1\textwidth]{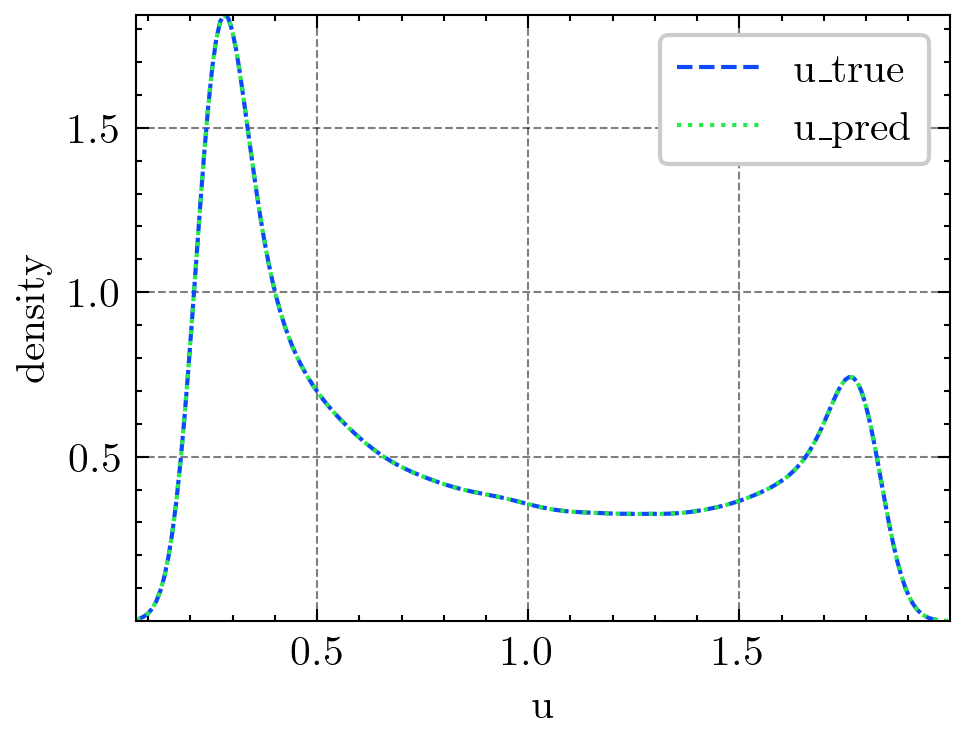}
        \label{fig:fk_density_10}
    \end{minipage}
    \caption{For equation of \eqref{eq:fokker_planck_equation}, comparison of the densities of eigenfunctions learned by IPMNN  and  the exact eigenfunctions  in (Top:) $d = 1$ (left) and $d= 2$ (right); (Bottom:) $d=5$ (left)  and $d=10$ (right). }
    \label{fig:fk_density}
\end{figure}

In this subsection we consider the linear Fokker-Planck equation with the periodic boundary condition in $\Omega=[0, 2\pi]^d$,
\begin{equation}
    \label{eq:fokker_planck_equation}
    \begin{array}{r@{}l}
        \begin{aligned}
             & -\Delta u - \nabla V \cdot \nabla u - \Delta Vu  =\lambda u, &  & \mbox{in} \enspace \Omega, \\
        \end{aligned}
    \end{array}
\end{equation}
where $V(\boldsymbol{x})$ is a potential function. The smallest eigenvalue is $\lambda=0$ and the corresponding eigenfunction is $u(\boldsymbol{x})=e^{-V(\boldsymbol{x})}$.
We choose $V(\boldsymbol{x})= \sin(\sum_{i=1}^dc_i \cos(x_i))$, where $c_i$ takes values in $[0.1, 1]$. 
The IPMNN is utilized to solve this problem with shifted operator due to the specificity of the smallest eigenvalue being 0. 
We choose $\alpha = 1$ and the parameters which are used to train IPMNN for the Fokker-Planck equation \eqref{eq:fokker_planck_equation} in different dimensions are summarized in  Table \ref{tab:params_ipmnn_fk}. 
To enforce the periodic boundary conditions as discussed above, we need to change the number of neurons in the input layer as $2dk$. In Figure \ref{fig:fk_line}, the smallest eigenfunction obtained by IPMNN, the exact eigenfunction  and the variation of loss are presented. The eigenfunction learned by IPMNN in two dimension, the exact solution and the absolute error between the NN solution and the true solution are shown in Figure \ref{fig:ipmnn_2D_heatmap3_fk}. It is obvious to see that IPMNN perfectly learns the eigenfunction in two dimension. We also implement our method in higher dimensions and we compare our results with that obtained by FBSDE \citep{han2020solving} in Table \ref{tab:error_ipmnn_fk}. The difference between exact eigenvalue and approximate eigenvalue is small enough to demonstrate the accuracy of IPMNN in different dimensions. 
The densities of eigenfunctions of equation \eqref{eq:fokker_planck_equation} in $d= 1$, $d = 2$, $d = 5$ and $d =10$ are shown in Figure \ref{fig:fk_density}. It is obvious to see that the densities of eigenfunctions learned by IPMNN perfectly fit the densities of exact eigenfunctions in all cases.


\begin{table}[htp]
\begin{center}
    \caption{Parameter settings of training IPMNN for Fokker-Planck equation \eqref{eq:fokker_planck_equation} in different dimensions.}
    \begin{tabular}{llllll}
        \hline\noalign{\smallskip}
        $d$ & $\alpha$ & $k$ & $N$    & $N_{epoch}$ & Layers of MLP           \\
        \hline
          1   & 1        & 3   & 10000  & 50000       & [6, 20, 20, 20, 20, 1]  \\
        2   & 1        & 3   & 20000  & 50000       & [12, 40, 40, 40, 40, 1] \\
         5   & 1        & 3   & 50000  & 50000       & [30, 60, 60, 60, 60, 1] \\
         10  & 1        & 3   & 100000 & 100000      & [60, 80, 80, 80, 80, 1] \\
        \hline
    \end{tabular}
    \label{tab:params_ipmnn_fk}
\end{center}
\end{table}

\begin{table}[htp]
\begin{center}
    \caption{Approximate eigenvalues for Fokker-Planck equation \eqref{eq:fokker_planck_equation} in different dimensions when exact eigenvalue $\lambda = 0$.}
    \begin{tabular}{lll}
        \hline\noalign{\smallskip}
        Method                      & $d$  & Approximate $\lambda$ \\
        \hline
        IPMNN                       & 1                  & 1.5497E-06            \\
        IPMNN                       & 2                  & 4.4227E-05            \\
        IPMNN                       & 5                  & 3.2902E-05            \\
        FBSDE \citep{han2020solving} & 5                  & 3.08E-03              \\
        IPMNN                       & 10                 & 1.0347E-04            \\
        FBSDE \citep{han2020solving} & 10                 & 3.58E-03              \\
        \hline
    \end{tabular}
    \label{tab:error_ipmnn_fk}
\end{center}
\end{table}


\subsubsection{Interior Eigenvalues Solved by IPMNN}
The method IPMNN is also able to find the interior eigenvalues and the associated eigenfunctions when some mild prior knowledge of the eigenvalue is provided. 
When any  $\alpha$ in equation \eqref{eq:Lu_equation1} is given,  
we solve for the nearest eigenvalue to $\alpha$ by using IPMNN.

\begin{figure}[htp]
    \centering
    \includegraphics[width=1\textwidth]{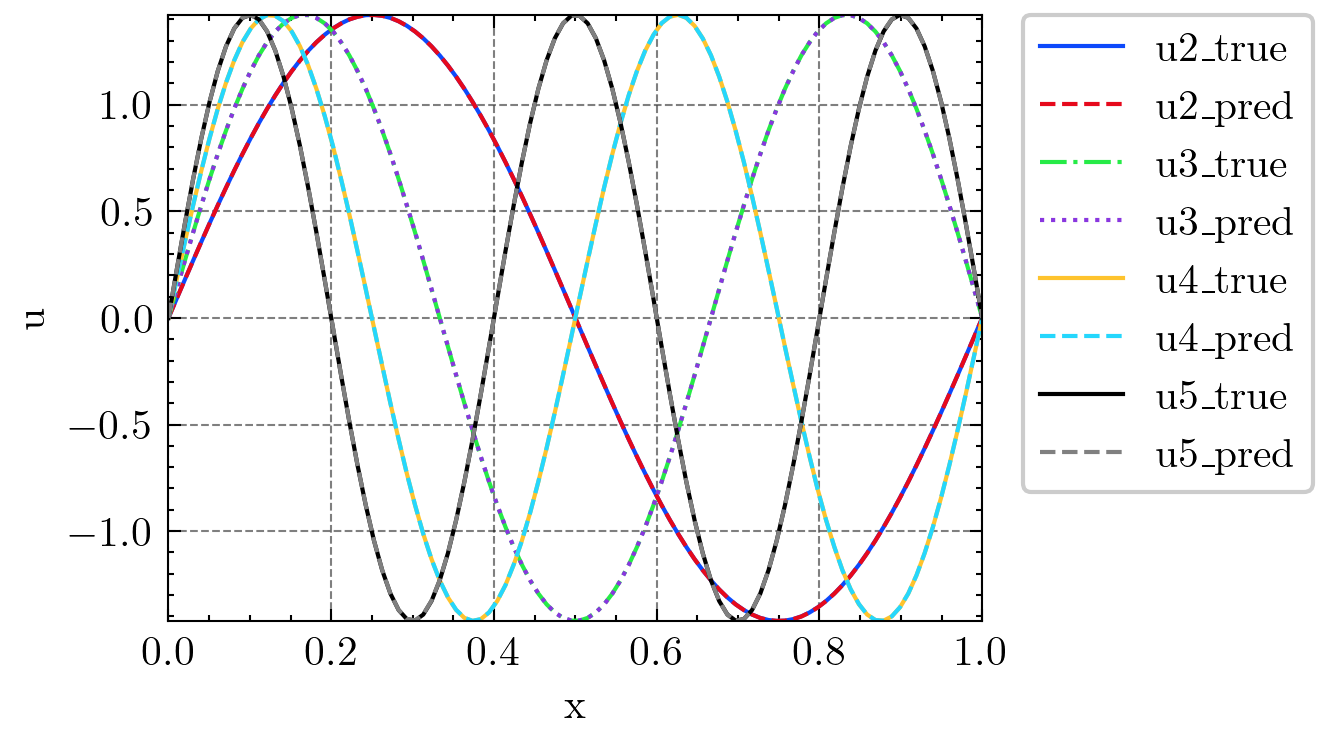}
    \caption{The eigenfunctions solved by IPMNN for different eigenvalues.}
    \label{fig:ipm_eigens_line}
\end{figure}

For simplicity, we consider the equation \eqref{eq:smallest_problem1} in one dimension. 
For given $\alpha$, the eigenvalues learned by IPMNN and the relative error are shown in Table \ref{tab:other_eigenvalues}. The comparison of eigenfunctions learned by IPMNN and the exact eigenfunctions is shown in Figure \ref{fig:ipm_eigens_line}. It is observed that the eigenvalues and eigenfunctions are both accurately solved by IPMNN.

\begin{remark}
    IPMNN is able to solve any eigenvalues when the prior knowledge $\alpha$ is given. However, if $\alpha$ is far from the eigenvalue that we want to find, it will often fail to learn it. Therefore, it is important to know the distribution of eigenvalues if we want to solve for some specified eigenvalues.
\end{remark}

\begin{table}[htp]
\begin{center}
    \caption{Exact eigenvalues and approximate eigenvalues learned by IPMNN for the equation \eqref{eq:smallest_problem1} with different $\alpha$.}
    \begin{tabular}{lllll}
        \hline\noalign{\smallskip}

                 $\alpha$ & Exact $\lambda$ & Approximate $\lambda$  & Relative error \\
        \hline
         36       & 39.4784         & 39.4853                    & 1.7334E-04     \\
         81       & 88.8264         & 88.8237                    & 3.0531E-05     \\
         144      & 157.9137        & 157.9074                    & 3.9804E-05     \\
         225      & 246.7401        & 246.7468                   & 2.7096E-05     \\        \hline
    \end{tabular}
    \label{tab:other_eigenvalues}
\end{center}
\end{table}

\section{Conclusions}
\label{sec:conclusions}
In this work, we propose PMNN  and IPMNN  which are neural networks combined with the power method and the inverse power method to solve linear eigenvalue problems. 
In similar spirits of traditional methods, our methods are presented for linear differential operators, the boundary conditions are enforced and the eigenvalue and eigenfunction are approximated iteratively.
Different from conventional numerical methods, we use AD to represent the differential operator. The eigenfunction is discretized to an eigenvector in the conventional methods, but we use the neural network to learn the eigenfunction. 
Another feature of our methods is the loss function which is different from the definitions in the usual form and is inspired from the power method and the inverse power method.  

Numerical experiments are carried out to verify the applicability and accuracy of the proposed methods for eigenvalue problems including high dimensional cases. We compare the predicted eigenfunctions with the exact solution in one and two dimensional problems. For higher dimensional problems, we compare the densities of eigenfunctions learned by our methods with the densities of the exact eigenfunctions. All the numerical results demonstrate that our methods can accurately learn the corresponding eigenfunctions. Additionally, we show the approximate eigenvalues for $d = 1, 2, 5$ and $10$. The eigenvalues are compared with the exact solutions, and also compared with the results obtained by other methods in the literature, like DRM and FBSDE. We get better results using our methods. Finally, 
the method IPMNN can also solve for interior eigenvalues with any given prior knowledge $\alpha$. 

Although good results are obtained, there are still some aspects need to be studied, which will be our future work. At first, we only discuss linear differential operators in this work. We would like to generalize our methods to nonlinear operators. 
Then, it will be a meaningful work for further investigate PINN for solving complex eigenvalue problems in a wide applications, such as in nuclear reactor physics domain, elastic acoustic problem, elastic viscoelastic composite structures, etc. And it will be interesting to investigate the efficiency of the proposed method for more realistic problems with much more complex geometry or materials etc.
Furthermore, if the dominant eigenvalue or the smallest eigenvalue has multiple values, our methods may fail. This is one direction that we will work on. 
Finally, the sampling algorithm and the number of sampling points are also important for our methods, which are very technical. 

\section*{Acknowledgments}
This research is supported part by the National Natural Science Foundation of China (No.11971020).
%
\section*{Data availability}
All data and codes used in this manuscript are publicly available at \url{https://github.com/SummerLoveRain/PMNN_IPMNN}.

 
 






\bibliographystyle{elsarticle-num}
\bibliography{ecrc-template}

\begin{thebibliography}{10}
\expandafter\ifx\csname url\endcsname\relax
  \def\url#1{\texttt{#1}}\fi
\expandafter\ifx\csname urlprefix\endcsname\relax\def\urlprefix{URL }\fi
\expandafter\ifx\csname href\endcsname\relax
  \def\href#1#2{#2} \def\path#1{#1}\fi

\bibitem{yu2018deep}
B.~Yu, W.~E, The deep {R}itz method: a deep learning-based numerical algorithm
  for solving variational problems, Communications in Mathematics and
  Statistics 6~(1) (2018) 1--12.

\bibitem{PINN}
M.~Raissi, P.~Perdikaris, G.~E. Karniadakis, Physics-informed neural networks:
  A deep learning framework for solving forward and inverse problems involving
  nonlinear partial differential equations, Journal of Computational physics
  378 (2019) 686--707.

\bibitem{lagaris1998artificial}
I.~E. Lagaris, A.~Likas, D.~I. Fotiadis, Artificial neural networks for solving
  ordinary and partial differential equations, IEEE transactions on neural
  networks 9~(5) (1998) 987--1000.

\bibitem{lee1990neural}
H.~Lee, I.~S. Kang, Neural algorithm for solving differential equations,
  Journal of Computational Physics 91~(1) (1990) 110--131.

\bibitem{van1995neural}
B.~P. Van~Milligen, V.~Tribaldos, J.~A. Jim{\'e}nez, Neural network
  differential equation and plasma equilibrium solver, Physical review letters
  75~(20) (1995) 3594.

\bibitem{cai2021physics}
S.~Cai, Z.~Wang, S.~Wang, P.~Perdikaris, G.~E. Karniadakis, Physics-informed
  neural networks for heat transfer problems, Journal of Heat Transfer 143~(6).

\bibitem{bai2022application}
Y.~Bai, T.~Chaolu, S.~Bilige, The application of improved physics-informed
  neural network (ipinn) method in finance, Nonlinear Dynamics 107~(4) (2022)
  3655--3667.

\bibitem{gao2022wasserstein}
Y.~Gao, M.~K. Ng, Wasserstein generative adversarial uncertainty quantification
  in physics-informed neural networks, Journal of Computational Physics 463
  (2022) 111270.

\bibitem{oszkinat2022uncertainty}
C.~Oszkinat, S.~E. Luczak, I.~Rosen, Uncertainty quantification in estimating
  blood alcohol concentration from transdermal alcohol level with
  physics-informed neural networks, IEEE Transactions on Neural Networks and
  Learning Systems (2022) 1--8\href
  {http://dx.doi.org/10.1109/TNNLS.2022.3140726}
  {\path{doi:10.1109/TNNLS.2022.3140726}}.

\bibitem{yang2019adversarial}
Y.~Yang, P.~Perdikaris, Adversarial uncertainty quantification in
  physics-informed neural networks, Journal of Computational Physics 394 (2019)
  136--152.

\bibitem{chen2020physics}
Y.~Chen, L.~Lu, G.~E. Karniadakis, L.~Dal~Negro, Physics-informed neural
  networks for inverse problems in nano-optics and metamaterials, Optics
  express 28~(8) (2020) 11618--11633.

\bibitem{kadeethum2020physics}
T.~Kadeethum, T.~M. J{\o}rgensen, H.~M. Nick, Physics-informed neural networks
  for solving inverse problems of nonlinear biot's equations: Batch training,
  in: 54th US Rock Mechanics/Geomechanics Symposium, OnePetro, 2020.

\bibitem{elhareef2022physics}
M.~H. Elhareef, Z.~Wu, Physics-informed neural network method and application
  to nuclear reactor calculations: A pilot study, Nuclear Science and
  Engineering (2022) 1--22.

\bibitem{buchan2013pod}
A.~Buchan, C.~Pain, F.~Fang, I.~Navon, A pod reduced-order model for eigenvalue
  problems with application to reactor physics, International Journal for
  Numerical Methods in Engineering 95~(12) (2013) 1011--1032.

\bibitem{diao2022spectral}
H.~Diao, H.~Li, H.~Liu, J.~Tang, Spectral properties of an acoustic-elastic
  transmission eigenvalue problem with applications, arXiv preprint
  arXiv:2210.16617.

\bibitem{chen2000integral}
Q.~Chen, Y.~Chan, Integral finite element method for dynamical analysis of
  elastic--viscoelastic composite structures, Computers \& Structures 74~(1)
  (2000) 51--64.

\bibitem{Golub1996}
G.~H. Golub, C.~F. Van~Loan, Matrix computations, Baltimore and London: John
  Hopkins University Press.

\bibitem{ben2020solving}
I.~Ben-Shaul, L.~Bar, N.~Sochen, Solving the functional eigen-problem using
  neural networks, arXiv preprint arXiv:2007.10205.

\bibitem{ben2023deep}
I.~Ben-Shaul, L.~Bar, D.~Fishelov, N.~Sochen, Deep learning solution of the
  eigenvalue problem for differential operators, Neural Computation 35~(6)
  (2023) 1100--1134.

\bibitem{bookalain2020}
A.~Hébert,
  \href{libgen.li/file.php?md5=cf327362bb2cc3a12162c4787b6717d8}{Applied
  Reactor Physics}, 3rd Edition, Presses internationales Polytechnique, 2020.
\newline\urlprefix\url{libgen.li/file.php?md5=cf327362bb2cc3a12162c4787b6717d8}

\bibitem{WUzeyun2022}
M.~H. ELHAREEF, Z.~WU, Extension of the pinn diffusion model to k-eigenvalue
  problems, American Nuclear Society (2022) (2022) 15--20.

\bibitem{JAGTAP2020113028}
A.~D. Jagtap, E.~Kharazmi, G.~E. Karniadakis,
  \href{https://www.sciencedirect.com/science/article/pii/S0045782520302127}{Conservative
  physics-informed neural networks on discrete domains for conservation laws:
  Applications to forward and inverse problems}, Computer Methods in Applied
  Mechanics and Engineering 365 (2020) 113028.
\newblock \href {http://dx.doi.org/https://doi.org/10.1016/j.cma.2020.113028}
  {\path{doi:https://doi.org/10.1016/j.cma.2020.113028}}.
\newline\urlprefix\url{https://www.sciencedirect.com/science/article/pii/S0045782520302127}

\bibitem{wang2022surrogate}
J.~Wang, X.~Peng, Z.~Chen, B.~Zhou, Y.~Zhou, N.~Zhou, Surrogate modeling for
  neutron diffusion problems based on conservative physics-informed neural
  networks with boundary conditions enforcement, Annals of Nuclear Energy 176
  (2022) 109234.

\bibitem{YANG2023109656}
Y.~Yang, H.~Gong, S.~Zhang, Q.~Yang, Z.~Chen, Q.~He, Q.~Li,
  \href{https://www.sciencedirect.com/science/article/pii/S0306454922006867}{A
  data-enabled physics-informed neural network with comprehensive numerical
  study on solving neutron diffusion eigenvalue problems}, Annals of Nuclear
  Energy 183 (2023) 109656.
\newblock \href
  {http://dx.doi.org/https://doi.org/10.1016/j.anucene.2022.109656}
  {\path{doi:https://doi.org/10.1016/j.anucene.2022.109656}}.
\newline\urlprefix\url{https://www.sciencedirect.com/science/article/pii/S0306454922006867}

\bibitem{jin2022physics}
H.~Jin, M.~Mattheakis, P.~Protopapas, Physics-informed neural networks for
  quantum eigenvalue problems, in: 2022 International Joint Conference on
  Neural Networks (IJCNN), IEEE, 2022, pp. 1--8.

\bibitem{han2020solving}
J.~Han, J.~Lu, M.~Zhou, Solving high-dimensional eigenvalue problems using deep
  neural networks: A diffusion monte carlo like approach, Journal of
  Computational Physics 423 (2020) 109792.

\bibitem{baydin2018automatic}
A.~G. Baydin, B.~A. Pearlmutter, A.~A. Radul, J.~M. Siskind, Automatic
  differentiation in machine learning: a survey, Journal of Marchine Learning
  Research 18 (2018) 1--43.

\bibitem{truhlar1972finite}
D.~G. Truhlar, Finite difference boundary value method for solving
  one-dimensional eigenvalue equations, Journal of Computational Physics 10~(1)
  (1972) 123--132.

\bibitem{simos1997finite}
T.~Simos, P.~Williams, A finite-difference method for the numerical solution of
  the schr{\"o}dinger equation, Journal of Computational and Applied
  Mathematics 79~(2) (1997) 189--205.

\bibitem{ishihara1977convergence}
K.~Ishihara, Convergence of the finite element method applied to the eigenvalue
  problem $\delta$u+ $\lambda$u= 0, Publications of the Research Institute for
  Mathematical Sciences 13~(1) (1977) 47--60.

\bibitem{ishihara1978mixed}
K.~Ishihara, A mixed finite element method for the biharmonic eigenvalue
  problems of plate bending, Publications of the Research Institute for
  Mathematical Sciences 14~(2) (1978) 399--414.

\bibitem{liang2001finite}
S.~Liang, X.~Ma, A.~Zhou, Finite volume methods for eigenvalue problems, BIT
  Numerical Mathematics 41~(2) (2001) 345--363.

\bibitem{dai2011finite}
X.~Dai, X.~Gong, Z.~Yang, D.~Zhang, A.~Zhou, Finite volume discretizations for
  eigenvalue problems with applications to electronic structure calculations,
  Multiscale Modeling \& Simulation 9~(1) (2011) 208--240.

\bibitem{talbot2005application}
C.~J. Talbot, A.~Crampton, Application of the pseudo-spectral method to 2d
  eigenvalue problems in elasticity, Numerical Algorithms 38~(1) (2005)
  95--110.

\bibitem{atkinson2009spectral}
K.~Atkinson, O.~Hansen, A spectral method for the eigenvalue problem for
  elliptic equations, arXiv preprint arXiv:0909.3607.

\bibitem{DW2002}
D.~Kincaid, D.~R. Kincaid, E.~W. Cheney, Numerical analysis: Mathematics of
  scientific computing, Americal Mathematical Society.

\bibitem{evans2010partial}
L.~C. Evans, Partial differential equations, American Mathematical Soc. 19.

\bibitem{berg2018}
J.~Berg, K.~Nystr$\ddot{o}$m, A unified deep artificial neural network approach
  to partial differential equations in complex geometries, Neuro computing 317
  (2018) 28--41.

\bibitem{lyu2020enforcing}
L.~Lyu, K.~Wu, R.~Du, J.~Chen, Enforcing exact boundary and initial conditions
  in the deep mixed residual method, arXiv preprint arXiv:2008.01491.

\bibitem{dong2021method}
S.~Dong, N.~Ni, A method for representing periodic functions and enforcing
  exactly periodic boundary conditions with deep neural networks, Journal of
  Computational Physics 435 (2021) 110242.

\bibitem{loh1996latin}
W.~L. Loh, On latin hypercube sampling, The annals of statistics 24~(5) (1996)
  2058--2080.

\end{thebibliography}







\end{document}